\documentclass[hidelinks,onefignum,final,onetabnum]{siamart250211}

\usepackage{lipsum}
\usepackage{amsfonts}
\usepackage{amssymb}
\usepackage{graphicx}
\usepackage{epstopdf}
\usepackage{pgfplots}
\usepackage{pgfplotstable}
\usepackage{amsmath}

\usepackage{algorithmic}
\ifpdf
  \DeclareGraphicsExtensions{.eps,.pdf,.png,.jpg}
\else
  \DeclareGraphicsExtensions{.eps}
\fi

\newsiamremark{remark}{Remark}
\newsiamremark{hypothesis}{Hypothesis}
\crefname{hypothesis}{Hypothesis}{Hypotheses}
\newsiamthm{claim}{Claim}
\newsiamthm{property}{Property}
\newsiamthm{problem}{Problem}
\newsiamremark{fact}{Fact}
\crefname{fact}{Fact}{Facts}

\headers{Krylov and Core Methods for MOP Recurrences}{A. Faghih, M. Rinelli, M. Van Barel, R. Vandebril and R. Vermeiren}

\title{Krylov and core transformation algorithms for an inverse eigenvalue problem to compute recurrences of multiple orthogonal polynomials}

\author{
Amin Faghih\thanks{Numerical Analysis and Applied Mathematics (NUMA) unit, Department of Computer Science, KU Leuven, Leuven, Belgium. \\
Amin Faghih: \email{amin.faghih@kuleuven.be}, 
Michele Rinelli: \email{michele.rinelli@kuleuven.be}, 
Marc Van Barel: \email{marc.vanbarel@kuleuven.be}, 
Raf Vandebril: \email{raf.vandebril@kuleuven.be}, 
Robbe Vermeiren: \email{robbe.vermeiren@kuleuven.be}}
\and Michele Rinelli\footnotemark[1]
\and Marc Van Barel\footnotemark[1]
\and Raf Vandebril\footnotemark[1]
\and Robbe Vermeiren\footnotemark[1] \textsuperscript{,}\thanks{Corresponding author}}

\usepackage{amsopn}
\DeclareMathOperator{\diag}{diag}

\ifpdf
\hypersetup{
  pdftitle={Inverse Eigenvalue Problem linked to Multiple Orthogonal Polynomials},
  pdfauthor={A. Faghih, M. Rinelli, M. Van Barel, R. Vandebril, R. Vermeiren}
}
\fi

\usepackage{subcaption}
\usepackage{bm}
\usepackage{tikz}
\usepackage{xcolor}
\usepackage{cite}
\usepackage{enumitem}

\hypersetup{
	colorlinks=true,
	linkcolor=blue,      
}
\usetikzlibrary{arrows.meta,calc,
	tikzmark,pgfplots.groupplots,external}
\pgfplotsset{compat=1.18}

\newcommand{\R}{\mathbb{R}} 
\renewcommand{\v}{\bm v} 
\newcommand{\w}{\bm w} 

\newcommand{\kryl}{\mathcal{K}} 
\newcommand{\dimk}{m} 

\newcommand\Tstrut{\rule{0pt}{2.6ex}}

\DeclareCaptionLabelFormat{opening}{#2}
\Crefname{ALC@unique}{Line}{Lines} 
\newcommand{\thicktimes}%
{\ooalign{$\bm\times$\cr$\bm\times$\cr\kern0.1pt$\bm\times$}}

\newcommand{\rev}[1]{#1} 

\begin{document}

\maketitle
\begin{abstract}

    In this paper, we develop algorithms for computing the recurrence coefficients corresponding to multiple orthogonal polynomials on the step-line. We reformulate the problem as an inverse eigenvalue problem, which can be solved using numerical linear algebra techniques. We consider two approaches: the first is based on the link with block Krylov subspaces and results in a biorthogonal Lanczos process with multiple starting vectors; the second consists of applying a sequence of Gaussian eliminations on a diagonal matrix to construct the banded Hessenberg matrix containing the recurrence coefficients. We analyze the accuracy and stability of the algorithms with numerical experiments on the ill-conditioned inverse eigenvalue problems have related to Kravchuk and Hahn polynomials, as well as on other better conditioned examples.
\end{abstract}



\section{Introduction}
\rev{Orthogonal polynomials form a fundamental tool in the analysis and development of numerical algorithms~\cite{Gautschi04,OlSlTo20,Liesenn}. We say that two polynomials $p,q$ are orthogonal with respect to a positive measure $\mu$ on $\mathbb{R}$ when $\int p(x)\,q(x)\,d\mu(x)= 0$. Classical families such as the Jacobi polynomials admit explicit formulas~\cite{Ism05}, but in many applications one deals with inner products that do not correspond to such families. In that case, the corresponding orthogonal polynomials must be computed numerically for which different methods exist~\cite{Gautschi04}. Several of these methods are based on numerical linear algebra techniques. For example, when the inner product is discrete, the recurrence coefficients corresponding to the orthogonal polynomials can be computed by solving an inverse eigenvalue problem~\cite{GraHar84}.
In this paper, we generalise some of these methods to the case of multiple orthogonal polynomials (MOPs).}

Recently, there has been a renewed interest in MOPs, which are an extension of the classical notion of orthogonal polynomials. This extension originates from simultaneous rational approximation, particularly from Hermite-Padé approximation, where multiple functions are approximated by rational functions with the same denominator~\cite{Apt98,AptMarRoc97,PAde}. MOPs are polynomials in one variable with orthogonality relations with respect to $r$ measures $\mu_{1},\mu_{2},\ldots,\mu_{r}$ \cite{Apt98,AptMarRoc97,Ceniceros2008,Ism05,VanaCou01}. In this text, we deal with
    positive discrete measures on $\mathbb{R}$~\cite{ArvCousVanAssche2003}, which also arise in quadrature rules for continuous measures~\cite{Vana24}. 
    
   Let $\vec{n}=(n_1, n_2, \dots, n_r) \in \mathbb{N}^r$ be a multi-index with size $ |\vec{n}| = n_1 + n_2 + \dots + n_r$. There exist two types of multiple orthogonal polynomials. Type I MOPs are tuples of polynomials $ (A_{\vec{n},1}, \dots, A_{\vec{n},r})$, where  $\deg A_{\vec{n},j} \leq n_j - 1$ for $j=1,\dots,r$, that satisfy the orthogonality conditions
    \begin{equation}
    \label{eq:type1}
        \int_\mathbb{R} \sum_{j=1}^{r} A_{\vec{n},j}(x) x^k d\mu_j(x) = 0, \quad  0\le k \le |\vec{n}| - 2.
    \end{equation}
    A type II multiple orthogonal polynomial $ P_{\vec{n}}$ is a polynomial of degree at most $|\vec{n}|$ that satisfies
    \begin{equation}
    \label{eq:type2}
        \int_\mathbb{R} P_{\vec{n}}(x) x^k d\mu_j = 0, \quad 0\le k \le n_j - 1, \quad 1 \le j \le r.
    \end{equation}
    
    Multiple orthogonal polynomials appear in various applications, such as rational approximation in number theory~\cite{Apt98}, evaluation of matrix functions~\cite{AlqRei18}, spectral theory of non-symmetric operators~\cite{AptKal98,Apt14}, special function theory~\cite{AptMarRoc97,VanaCou01}, random matrix theory~\cite{Kuijlaars2010}, eigenvalues of products of random matrices, non-intersecting Brownian motions, and combinatorial problems linked to random tilings~\cite{Ism05,Kuijlaars2010}.

In 1969, Golub and Welsch~\cite{GoWe69} illustrated that in Gaussian quadrature (i.e., the case $r=1$), the nodes $z_i$ and weights $w_i$ can be expressed in terms of eigenvalues and first components of the eigenvectors of a Jacobi matrix, respectively. \rev{This can be explained by looking at the orthogonal polynomials $\{p_k\}_{k=0}^{N-1}$ corresponding to the discrete inner product $\langle p,q\rangle_N = \sum_{i=1}^N |w_i|^2 p(z_i)q(z_i)$. Enforcing orthogonality  and using the recurrence relations naturally leads to the following relation~\cite{RafMarcNicola08}
\begin{equation*}
    ZQ=QT,
\end{equation*}
with $Z$ a diagonal matrix containing the nodes, $Q \in \mathbb{R}^{N{\times} N}$ an orthogonal matrix containing the evaluations of $p_k$ in the nodes, and $T\in \mathbb{R}^{N{\times} N}$ a Jacobi matrix containing the coefficients appearing in the three-term recurrence relation.}
    MOPs satisfy several recurrence relations depending on the path chosen for the multi-indices in $\mathbb{N}^r$. In this paper, we focus on the so-called step-line recurrence \rev{(see \cref{Sec3})}, an $(r+2)$-term relation, resulting in a banded upper Hessenberg matrix $H_N$ taking the role of the Jacobi matrix~\cite{FilHanVanass15,Ism05,VanaCou01}.\footnote{In the literature, the recurrence relations are commonly written with a banded, lower Hessenberg matrix. However, we choose to work with the transpose of this matrix, as this matches more closely to the literature in numerical linear algebra, on which we will heavily rely.}
    As a consequence, the type I and type II MOPs are linked to the right and left eigenvectors of the recurrence matrix $H_N$, respectively; see, e.g.,~\cite{CouVanass05}. 
    We propose numerical algorithms to compute
    these recurrence relations and recover the corresponding MOPs.  
    
    The Golub-Welsch approach was extended to multiple Gaussian quadrature by Borges~\cite{Bor94}, Coussement et al.~\cite{CouVanass05}, Van Assche~\cite{Vana24}, and Laudadio et al.~\cite{LauMasVand24}, where the goal is to compute optimal nodes $\{z_i\}_{i=1}^N$ and weights $\{\alpha_{j,i}\}_{i=1}^N$ for $j=1,2,\ldots,r$, to approximate integrals of a function $f$ with respect to the different measures:
\begin{equation*}
    \int_{\mathbb{R}} f(x)d\mu_{j}(x) \approx \sum_{i=1}^N \alpha_{j,i} f(z_i).
\end{equation*} 
We will be concerned with the opposite problem: given the nodes and weights of some discrete measures, we desire to retrieve the corresponding recurrence matrix. As explained in the previous paragraph, the nodes and weights are linked to the eigenvalues and eigenvectors of this recurrence matrix. Hence, we have an instance of an IEP. \cite{Chu98,BolGol87}. 

We present two main algorithms for solving the IEP associated with MOPs. The first exploits the connection between MOPs and the pair of biorthogonal bases generated by a Lanczos-type algorithm with multiple starting vectors, which is a specialized version of the method described by Aliaga et al.~\cite{AliBolFreHer00}. 

The second is a core transformation algorithm which enforces the conditions specified in the IEP by iteratively executing Gaussian transformations on the diagonal matrix of nodes, generalizing the techniques developed by Mach, Van Barel, and Vandebril~\cite{MachVanbVand14}. Both methods have been proven particularly effective in similar problems: 
generating orthogonal polynomials on the real line by solving an IEP associated with a Jacobi matrix~\cite{GraHar84,Rut63}, constructing a set of biorthogonal polynomials following from the solution of a tridiagonal IEP \cite{VB21}, generating a sequence of (bi)orthogonal rational functions by solving a (tridiagonal) Hessenberg pencil IEP~\cite{VanbVanbVand22}, generating Sobolev orthogonal polynomials involving a Hessenberg IEP~\cite{Vanb23}, and generating Sobolev orthogonal rational functions through a Hessenberg pencil IEP \cite{Faghih2025Sobolev,FaVBVBVa24}. Furthermore, if one aims to do updating~\cite{VanbVanbVand22} and downdating~\cite{VanBVanBVand24}, i.e. adding or removing a node in the discrete inner product, the core transformation approach is preferred. A third method exists, based on the Cholesky factorization of the moment matrix~\cite{VB21}. We briefly outline this approach, but do not study it in detail because of its less favorable numerical properties. 
\rev{For classical orthogonal polynomials, Laurie~\cite{Laurie99} devised another strategy based on the connection between classical orthogonal polynomials and continued fractions. We are not aware of any extension of Laurie’s approach to MOPs. Possible developments are discussed, with further references, in \cref{sec:conclusion}, but are beyond the scope of the paper.}

The computation of the step-line recurrence matrix $H_N$ was already investigated by Filipuk, Haneczok and Van Assche~\cite{FilHanVanass15}, who proposed an algorithm to transform the recurrence coefficients of the individual measures into those for the multiple orthogonal polynomials (MOPs). Additionally, Milovanović and Stanić~\cite{MilStan03} introduced a discretized Stieltjes-Gautschi procedure for computing this recurrence matrix. This procedure expresses the entries of $H_N$ in terms of the different inner products, which are then discretized using the corresponding Gaussian quadrature rules.\footnote{Every measure $\mu_j$ corresponds to an inner product $\left\langle f,g \right\rangle_{\mu_j} = \int_{\mathbb{R}}f(x)g(x)d\mu_j(x)$ of the functions $f$ and $g$ \rev{in an appropriate function space}.} Note that in our setting, this simplifies due to the assumption of discrete measures. Neither of these studies, however, explored the connection with Krylov subspaces or analyzed the numerical accuracy of the resulting methods. A systematic comparison of the performance of all these algorithms would be interesting for further research.

The paper is organized as follows: in \cref{Sec2}, the uniqueness conditions of type I and type II MOPs in a system of discrete measures are discussed. \Cref{Sec3} discusses the recurrence relations and formulates this as a matrix relation that leads to the formulation of the IEP. 
The link between Type I and block Krylov subspaces, and between Type II MOPs and standard Krylov subspaces, is explored in \cref{sec:Krylov}. Also, we provide a Lanczos-type algorithm and discuss breakdown and reorthogonalization. In \cref{Sec5}, we describe the core transformation algorithm and illustrate its functioning with a small example. The outcome of both algorithms is the structured matrix whose elements are the recurrence coefficients for generating a sequence of biorthogonal type I and type II MOPs. 
In \cref{sec:numericalexperiments}, we test the correctness and accuracy of the proposed algorithms on the IEPs induced by the Kravchuk and Hahn MOPs, for which analytic formulae are derived in~\cite{ArvCousVanAssche2003}. The poor performance is attributed to the ill-conditioning of the IEPs, as revealed by the Multiprecision toolbox Advanpix in MATLAB~\cite{advanpix2015}; on the other hand, we get favorable results on artificial inner products with randomly generated weights, which turn out to be better conditioned. All MATLAB codes required to reproduce the experiments are publicly available.\footnote{\url{https://gitlab.kuleuven.be/numa/public/IEP-MOP}} Finally, \cref{sec:conclusion} gives a summary and discusses possible directions for future research.

\section{Uniqueness of type I and type II MOPs}\label{Sec2}
In this section, we analyze the existence and uniqueness of solutions to the linear systems \cref{eq:type1} and \cref{eq:type2} by introducing an algebraic Chebyshev system (AT-system).
    
For the type \MakeUppercase{\romannumeral 2} case, the usual requirement is that every possible solution of \cref{eq:type2} has exactly degree $|\vec{n}|$. If this holds, we call $\vec{n}$ a \emph{normal index}. Since $P_{\vec{n}}$ is a polynomial of degree $|\vec{n}|$, equation \cref{eq:type2} defines a homogeneous linear system of $|\vec{n}|$ equations with $|\vec{n}|+1$ unknowns. The coefficient matrix $D^{(II)}$ of this system contains the moments of $\mu_j$, i.e., $m^{(j)}_{l}=\int_{\mathbb{R}} x^l d\mu_{j}$. Let us define the block matrix ${M}_{\vec{n}}$ of size $|\vec{n}| \times |\vec{n}|$ as $
{M}_{\vec{n}}= \begin{bmatrix}M^{(1)}_{n_{1}} & M^{(2)}_{n_{2}} & \ldots & M^{(r)}_{n_{r}}\end{bmatrix}$, with the rectangular blocks 
\begin{equation}\label{eq:momentmatrix}
    M^{(j)}_{n_{j}} = \begin{bmatrix}
m_0^{(j)}&m_1^{(j)} &\ldots& m_{n_{j}-1}^{(j)} \\[0.5em]
m_1^{(j)} & m_2^{(j)}&\ldots& m_{n_{j}}^{(j)} \\

\vdots & \vdots&&\vdots \\
m_{|\vec{n}|-1}^{(j)} &m_{|\vec{n}|}^{(j)}&\ldots& m_{|\vec{n}|+n_j-2}^{(j)}\\[0.5em]
\end{bmatrix} \in \mathbb{R}^{|\vec{n}| \times n_{j}}.
\end{equation}
Note that $M_{\vec n}^{T}$ equals the coefficient matrix $D^{(II)}$ with its last column removed. Thus, $\vec n$ is a normal index precisely when $D^{(II)}$ has full row rank \cite{CouVanass05,ArvCousVanAssche2003}.

A common requirement for a good definition of type I MOPs is that each $A_{\vec{n},j}$ has exactly degree $n_j - 1$. Note that the coefficient matrix $D^{(I)}$ of linear system~\cref{eq:type1} coincides with ${M}_{\vec{n}}$ having its last row removed. Thus, the above requirement is satisfied if $D^{(I)}$ has full row rank. 

To have a unique solution for systems \cref{eq:type1} and \cref{eq:type2}, we need one additional requirement. For type \MakeUppercase{\romannumeral 2}, we require the polynomials $P_{\vec{n}}$ to be monic. For type \MakeUppercase{\romannumeral 1}, the following normalization condition is imposed
\begin{equation}\label{eq:normalizationcondition}
        \int_\R \sum_{j=1}^{r} A_{\vec{n},j}(x)x^{|\vec{n}|-1} d\mu_j(x) = 1.
\end{equation} 
To find more details on the theory of MOPs, we refer the readers to the work of Coussement and Van Assche \cite{CouVanass05}, and Ismail \cite{Ism05}.

In this paper, we restrict ourselves to discrete positive measures on $\R$
\begin{equation}\label{MEASURE}
    \mu_j = \sum_{i=1}^{N} \alpha_{j,i} \delta_{z_{i}}, \; \; \alpha_{j,i} > 0, \; z_{i} \in \R, \; N \in \mathbb{N}, \; 1 \le j \le r,
\end{equation}
where all $z_{i}$ are distinct, and $\delta_{z_{i}}$ denotes the Dirac measure at $z_{i}$. The support of these discrete measures is defined as $\text{supp}(\mu_j) = \{z_{i}\}_{i=1}^{N}$. Also denote $\Delta \subset \mathbb{R}$ as an interval that contains $\operatorname{supp}(\mu_j)$ for all $j$. 
	\subsection{AT-system}\label{subsec:AT-system}
	An AT-system is a particular set of measures such that every multi-index $\vec{n}$ is normal. This subsection closely follows the discussion of Arves\'u, Coussement, and Van Assche~\cite{ArvCousVanAssche2003}.
	\begin{definition}[Chebyshev system]
		\sloppy A linearly independent system of $n$ basis functions $\alpha_1(x), \ldots, \alpha_n(x)$ on $\Delta \subset \R$ is a Chebyshev system if every \textit{non-trivial} linear combination $\sum_{k=1}^{n}a_k\alpha_k(x)$, with all $a_k \in \R$, has at most $n-1$ zeros in $\Delta$.
	\end{definition}
	\begin{definition}[AT-system \cite{ArvCousVanAssche2003}] \label{def:AT-system}
		An AT-system consists of $r$ positive discrete measures $\{\mu_j\}_{j=1}^r$ (cf. \cref{MEASURE}) with $\operatorname{supp}(\mu_j) = \{z_i\}_{i=1}^N$ for all $j$, such that there exist $r$ continuous functions $\alpha_j(x)$ on the interval $\Delta$ for which $\alpha_j(z_i) = \alpha_{j,i}$ for $j=1,\ldots,r,\; i=1,\ldots,N$ and the following $|\vec{n}|$ functions 
	\begin{equation*}
		\begin{aligned}
			\alpha_1(x), x\alpha_1(x), &\ldots, x^{n_1-1}\alpha_1(x) \\
			\alpha_2(x), x\alpha_2(x), &\ldots, x^{n_2-1}\alpha_2(x) \\
			&\vdots \\
			\alpha_r(x), x\alpha_r(x), &\ldots, x^{n_r-1}\alpha_r(x) 
		\end{aligned}
	\end{equation*} 
	give a Chebyshev system for each $\vec{n}$ where $|\vec{n}| \leq N$.

	\end{definition} 
    
	The condition $\lvert \vec{n} \rvert\leq N$ is necessary for ensuring uniqueness. Indeed, suppose that $m = \lvert \vec{n} \rvert - N  > 0$. Then we can take the polynomial $P_{\vec{n}}(x) = \prod_{i=1}^{N} (x-z_{i})\prod_{i=1}^{m} (x-a_i)$ which will be a solution of system \cref{eq:type2} for all $a_1,\ldots,a_m$. It is important to note that if $|\vec{n}| = N$, then $P_{\vec{n}}(x) = (x-z_1)(x-z_2)\cdots(x-z_N)$.

    The following theorem is a direct consequence of \cite[Theorem~2.1]{ArvCousVanAssche2003}, and ensures the uniqueness of MOPs in an AT-system.
	\begin{theorem}[Uniqueness]\label{theorem:AT}
		Assume we have an AT-system with $r$ positive discrete measures on $\Delta$. Then for $|\vec{n}| \leq N$
        \begin{itemize}
        \item The coefficients of the type I polynomials $(A_{\vec{n},1}, \dots, A_{\vec{n},r})$, $ \deg A_{\vec{n},j} = n_j - 1 $, as the solution of the linear system \cref{eq:type1} along with the normalization \eqref{eq:normalizationcondition}, can be uniquely determined.
        \item The linear system \cref{eq:type2} has a unique set of coefficients defining the type II monic orthogonal polynomial $ P_{\vec{n}}(x) $ of degree $|\vec{n}|$ as its solution.
        \end{itemize}
    \end{theorem}

	\section{Inverse eigenvalue problem}\label{Sec3}
    In this section, we discuss the connections between type I and type II MOPs, starting with their recurrence relations, followed by the biorthogonality property. Building upon these connections, we will subsequently formulate the problem of computing the recurrence coefficients as an inverse eigenvalue problem.

    \subsection{The step-line recurrence relation}\label{sec:steplinerecurrence}
    Multiple orthogonal polynomials satisfy a number of recurrence relations depending on how the indices $n_{j}$ are chosen. In this paper, we are concerned with the so-called step-line recurrence relation; see \cite[Sec.~23.1.4]{Ism05} for a discussion on other recurrence relations. Let us define the monic type II MOPs on the step-line by \begin{equation*}
        P_{n}(x)=P_{\vec{n}}(x), \quad \text{for} \quad n=kr+\ell,
    \end{equation*}
where $r$ is the number of measures, and $\vec{n}=(\underbrace{k+1, \ldots, k+1}_{\ell}, k, \ldots, k)$.
    It is known that the step-line monic polynomials $P_n(x)$ satisfy an $(r+2)$-term recurrence relation; see e.g., \cite{ArvCousVanAssche2003,CouVanass05,FilHanVanass15,Ism05,Vana24,VanaCou01},
    \begin{equation}\label{eq:recurrencetype2}
    xP_n(x)=P_{n+1}(x)+\sum_{j=0}^{r}a_{n,j}P_{n-j}(x), \quad \quad 0 \le n \le N-1,
    \end{equation}
    with initial conditions $P_{0} \equiv 1$ and $P_{j} \equiv 0$ for $j=-1,-2,\ldots,-r$.
    
    Let us define the discrete measure $\mu(x) = \sum_{i=1}^N \delta_{z_i}$, so that $d\mu_j(x) = \alpha_j(x) d\mu(x)$ with $\alpha_j(z_i) = \alpha_{j,i}$ and zero elsewhere.\footnote{$\mu_j$ is absolutely continuous with respect to $\mu$.} We define the type I functions as
    \begin{equation}\label{eq:type1functions}
        Q_{\vec{n}}(x) = \sum_{j=1}^{r} A_{\vec{n},j}(x) \alpha_j(x).
    \end{equation}
    Note that, with this notation, the orthogonality relations \cref{eq:type1} and the normalization condition \cref{eq:normalizationcondition} for type I MOPs become 
    \begin{equation} \label{NEWTYPE1}  
     \int_{\mathbb{R}} Q_{\vec{n}}(x) x^k  d\mu(x) = \sum_{i=1}^{N} Q_{\vec{n}}(z_{i}) z_{i}^k =\begin{cases}
    0,   & 0\leq k \leq |\vec{n}| - 2, \\
    1, & k = |\vec{n}| - 1.
    \end{cases}
    \end{equation}
    The type I functions~\cref{eq:type1functions} fulfill a recurrence relation similar to \cref{eq:recurrencetype2}, involving the same coefficients~\cite{Ism05,Vana24}:
    \begin{equation}\label{eq:recurrencetype1}
            xQ_n(x)=Q_{n-1}(x)+\sum_{j=0}^{r}a_{n+j-1,j}Q_{n+j}(x), \quad \quad 1\le n \le N-r,
    \end{equation}
    with initial conditions $Q_{0} \equiv 0$ and $Q_{1},Q_{2},\ldots,Q_{r}$.

    For clarity and ease of implementation, we restrict ourselves to the case $r=2$ for the remainder of the paper. 
    For $r=2$, the type I functions \cref{eq:type1functions} take the form 
    \begin{equation}\label{type1:r=2}
       Q_{n}(x)=\alpha_{1}(x)A_{n,1}(x)+\alpha_{2}(x)A_{n,2}(x),
    \end{equation}
    with $Q_n(x) := Q_{\vec{n}}(x)$, $A_{n,1}(x) := A_{\vec{n},1}(x)$ and $A_{n,2}(x) := A_{\vec{n},2}(x)$  where 
    \begin{align}
    \vec{n} = 
    \begin{cases}
    (k, k),   & n = 2k, \\
    (k+1, k), & n = 2k + 1.
    \end{cases}
    \label{eq:step-line-index-r=2}
    \end{align}
    In this setting, the step-line recurrence relation \eqref{eq:recurrencetype2} becomes
    \begin{equation}
        \label{eq:recurrence-typeII-polynomials}
        x P_n(x) = P_{n+1}(x) + b_n P_n(x) + c_n P_{n-1}(x) + d_n P_{n-2}(x), \quad 0 \le n\le N-1,
    \end{equation}
    while \eqref{eq:recurrencetype1} is written as
    \begin{equation}
        x Q_n(x) = Q_{n-1}(x) + b_{n-1} Q_n(x) + c_n Q_{n+1}(x) + d_{n+1} Q_{n+2}(x),\quad 1 \le n\le N-r. \label{eq:recurrence-typeI-functions}
    \end{equation}
    The coefficients of the recurrence relation \cref{eq:recurrence-typeII-polynomials} corresponding to $\{P_{n}(x)\}_{n=0}^{N-1}$ can be stored in a banded upper Hessenberg matrix
   \begin{equation}\label{eq:Hessenberg-recurrence-matrix}
   H_{N}=
    \begin{bmatrix}
    b_0 & c_1 & d_2 & 0 & 0&0&\cdots & 0 \\
    1 & b_1 & c_2 & d_3 & 0 &0& \cdots & 0 \\
    0 & 1 & b_2 & c_3 & d_4 &0& \cdots & 0 \\
    0 & 0 & 1 & b_3 & c_4 & d_5 & \cdots & 0 \\
    \vdots &  & \ddots & \ddots & \ddots & \ddots & \ddots &  \\
    0 & \cdots & 0 & 0 & 1 & b_{N-3} & c_{N-2} & d_{N-1} \\
    0 & \cdots & 0 & 0 & 0 & 1 & b_{N-2} & c_{N-1} \\
    0 & \cdots & 0 & 0 & 0 & 0 & 1 & b_{N-1}
    \end{bmatrix},
    \end{equation}
    allowing us to write the step-line recurrence relation for the type II polynomials as
    \begin{align}\label{eq:type2recurrence} 
        x \left[ P_0(x)  \; \ldots \; P_{N-1}(x) \right] 
        + \left[ 0 \; \ldots \; P_N(x) \right]= \left[ P_0(x) \; \ldots \; P_{N-1}(x) \right] H_N .
    \end{align}

   Recall that $P_N(x)$ is the unique monic polynomial of degree $N$ whose roots are the nodes $\{z_{i}\}_{i=1}^N$ associated with the measures $\mu_j$, cf.~\eqref{MEASURE}. Thus, evaluating the above relation in $z_{i}$ gives the following property.
    \begin{property} \label{property:lefteigenvector}
          If $P_{N}(z_{i})=0$, then $z_{i}$ is an eigenvalue of $H_{N}$ corresponding to left eigenvector $\bm{v}_{i,N}=\left[P_0(z_{i}) \; P_1(z_{i}) \; \ldots \; P_{N-1}(z_{i})\right] $.
    \end{property}
    Following this property, the roots of $P_n$, $n<N$ coincide with the eigenvalues of the leading $n\times n$ principal submatrix of $H_N$. These roots satisfy particular interlacing properties (see \cite{Interlacing}), but are not, in general, the eigenvalues of the full matrix $H_N$.
    
    By evaluating equations \cref{eq:type2recurrence} in the nodes $z_i$, we obtain
    \begin{equation} \label{eq:recurrence-V-matrix}
       ZV= V H_N, 
    \end{equation}
    with 
    \begin{align}\label{eq:VZ}
        \begin{array}{cc}
            V =
            \begin{bmatrix}
                P_0(z_1) & P_1(z_1) & \cdots & P_{N-1}(z_1) \\
                P_0(z_2) & P_1(z_2) & \cdots & P_{N-1}(z_2) \\
                \vdots   & \vdots   &        & \vdots       \\
                P_0(z_N) & P_1(z_N) & \cdots & P_{N-1}(z_N)
            \end{bmatrix},\; \;
            Z = 
            \begin{bmatrix}
                z_1 & & &\\
                 & z_2&& \\
                    & &\ddots &\\
                   & & & z_N
            \end{bmatrix}.
        \end{array}
    \end{align}

    \begin{remark}
    \label{rem:nonzero-second-superdiagonal}
        The type I MOPs satisfy the same recurrence relation \cref{eq:recurrence-typeI-functions} of the type I functions (see, e.g., \cite[Sec.~23.1.4]{Ism05}):
        \begin{equation}
            \label{eq:recurrence-typeI-polynomials}
            xA_{n,j}(x) = A_{n-1,j}(x) + b_{n-1}A_{n,j}(x) + c_n A_{n+1,j}(x) + d_{n+1}A_{n+2,j}(x),
        \end{equation}
        for $n\geq 1$ and $j=1,2$. The fact that the measures form an AT-system implies that $d_{n+1} \ne 0$ for all $n = 1, \dots, N-2$. Indeed, if $n=2k$ is even, we get from \cref{eq:step-line-index-r=2} that $\deg A_{n,2} = \deg A_{n+1,2} = k-1$ and $\deg A_{n+2,2}= k$. Hence, if we consider \cref{eq:recurrence-typeI-polynomials} for $j=2$, the left-hand side is a polynomial of degree $k+1$, while the right-hand side has degree $k+1$ if and only if $d_{n+1}\neq 0$. If $n=2k+1$ is odd, the same conclusion follows by observing that $\deg A_{n,1} = \deg A_{n+1,1} = k$ and $\deg A_{n+2,1} = k+1$. 
    \end{remark}
    \subsection{The inverse eigenvalue problem}
In this section, we formalize
the problem of generating MOPs associated with a discrete inner product as an inverse eigenvalue problem. Before doing this, we require one additional result establishing a biorthogonality relation between the type I functions and the type II monic polynomials. Given the discrete positive measures $\mu_{j}$ as in \cref{MEASURE}, we define a discrete inner product for type I functions and type II monic polynomials on the step-line, with real, positive weights  $\{ \alpha_{j,i}\}_{i=1}^N $, $j=1,2$, and distinct nodes $ \{ z_i \}_{i=1}^N $, $ z_i \in \mathbb{R} $, as 
\begin{equation}\label{INNER}
\langle P_{n},Q_{m}\rangle _N=\sum_{i=1}^N P_n(z_i)  Q_{m}(z_i)  =\sum_{i=1}^N P_n(z_i)  \big(\alpha_{1,i}A_{m,1}(z_i)+\alpha_{2,i}A_{m,2}(z_{i})\big) .
\end{equation}
    \begin{proposition}[Chapter 23 of \cite{Ism05}]\label{prop:bior}
        Given the type I functions as defined in \cref{type1:r=2} and the type II monic polynomials on the step-line, the following biorthogonality with respect to $\langle .,.\rangle _N$ holds:
        \begin{equation}\label{BIORTHOGONALITY}
         \langle P_{n},Q_{m}\rangle _N= \left\{\begin{array}{l}
            0,\quad \quad m\leq n ,\\
            0,\quad \quad n\leq m -2,\\
            1,\quad \quad m=n+1.
            \end{array}\right.
        \end{equation}
    \end{proposition}
    This result follows directly from the orthogonality conditions of type I and type II polynomials. The case where $m=n+1$ holds because of the additional normalization relation for type I \cref{eq:normalizationcondition}. In order to reformulate \cref{prop:bior}, we contruct a matrix $W$ of size $N \times N$ as
    \begin{equation*}
    W =
            \begin{bmatrix}
                Q_1(z_1) & Q_2(z_1) & \cdots & Q_N(z_1) \\
                Q_1(z_2) & Q_2(z_2) & \cdots & Q_N(z_2) \\
                \vdots   & \vdots   &        & \vdots   \\
                Q_1(z_N) & Q_2(z_N) & \cdots & Q_N(z_N)
            \end{bmatrix}.
            \end{equation*}
    The relation \cref{BIORTHOGONALITY} now becomes 
    \begin{equation}\label{eq:biorth}
        W^TV = I_N,
    \end{equation} 
    with $I_{N}$ the identity matrix of size $N \times N$. By combining the biorthogonality relation \cref{eq:biorth} with equation \cref{eq:recurrence-V-matrix}, the identities $ZW= W H_{N}^{T}$ and $H_N = W^T Z V$ are both equivalent to relation \cref{eq:recurrence-V-matrix}. Moreover, \rev{the $i$th column of $W^T$} is a right eigenvector of $H_N$ corresponding to the eigenvalue $z_i$.  The first and the second components of these right eigenvectors can be retrieved using the following result from Van Assche~\cite[Theorem $2$]{Vana24}.\footnote{Van Assche used this result in the context of Gaussian quadrature. He assumed that the first and second components of the right eigenvectors were given and calculated the optimal quadrature weights $\alpha_{1,j}$, $\alpha_{2,j}$.}
     \begin{theorem}\label{theorem:VanAsscheThm2}
     Given nodes $\{z_{i}\}_{i=1}^{N}$, and weight vectors $\bm{\alpha}_1 = \begin{bmatrix} \alpha_{1,1} \ldots \alpha_{1,N} \end{bmatrix}^T$ and $\bm{\alpha}_2 = \begin{bmatrix} \alpha_{2,1} \ldots \alpha_{2,N} \end{bmatrix}^T$ from the discrete measures $\mu_1$ and $\mu_2$ in equation \cref{MEASURE}, one can compute the first and second column of $W$ as 
    \begin{eqnarray*}
        \bm{\alpha}_{1}&=&d_{1}\bm{w}_{1}, \\
        \bm{\alpha}_{2}&=&d_{2}\bm{w}_{1}+d_{3}\bm{w}_{2},
    \end{eqnarray*}
    where $d_{1}$, $d_{2}$, and $d_{3}$ are constants given by 
      \begin{eqnarray*}
    d_1=\int_{\mathbb{R}}P_{0}(x)d\mu_1(x)=\sum_{i=1}^N \alpha_{1,i}P_{0}(z_i), \quad d_2=\int_{\mathbb{R}}P_{0}(x)d\mu_2(x)=\sum_{i=1}^N \alpha_{2,i}P_{0}(z_i), \\
   d_3=\int_{\mathbb{R}}P_{1}(x)d\mu_2(x)=\sum_{i=1}^N \alpha_{2,i}P_{1}(z_i).\hspace{+3cm}
   \end{eqnarray*}
      \end{theorem}

We can now state the main problem formally. We start with the general formulation.
    \begin{problem}[Generating MOPs]\label{Prob:wHrLS}
        Let a discrete inner product $\langle .,. \rangle _N$ as in \cref{INNER} be given, with nodes $\{z_{i}\}_{i=1}^{N}$ and a pair of weights $\{\alpha_{1,i}\}_{i=1}^{N}$, and $\{\alpha_{2,i}\}_{i=1}^{N}$. Compute a sequence of biorthogonal type I functions $\{Q_{n}\}_{n=1}^{N}=\{\alpha_{1}A_{n,1}+\alpha_{2}A_{n,2}\}_{n=1}^{N}$ and type II monic polynomials $\{P_{n}\}_{n=0}^{N-1}$ at the given nodes $z_{i}$ satisfying \cref{BIORTHOGONALITY}.
    \end{problem}
\cref{Prob:wHrLS} can be reformulated as an inverse eigenvalue problem. Here and throughout the rest of this paper, we denote the recurrence matrix by $H$ when the dependence of the size is clear from the context.
\begin{problem}\label{Prob:IEP}
Given a diagonal matrix $Z =
\diag(\{z_{i}\}_{i=1}^{N}) \in \mathbb{R}^{N\times N}$ with distinct nodes, and vectors $\bm{w}_{1},\bm{w}_{2}, \bm{v}_{1} \in \mathbb{R}^{N}$ satisfying
\begin{equation}\label{eq:initialdatarelations}
    \bm{w}_1^T \bm{v}_1 = 1, \quad \bm{w}_2^T \bm{v}_1 = 0, \quad \text{and} \quad \bm{w}_2^T Z \bm{v}_1 = 1,    
\end{equation}
construct \rev{an upper Hessenberg} matrix $H\in \mathbb{R}^{N \times N}$, as in \cref{eq:Hessenberg-recurrence-matrix}, and a pair of matrices $ W, V \in \mathbb{R}^{N \times N} $ such that
\begin{enumerate}[label=(D\arabic*)]
	\item\label{cond:1} The first two columns of $W$ equal $\bm{w}_1$ and $\bm{w}_2$ and the first column of $V$ is $\bm{v}_1$:
    \[
	\begin{bmatrix}
		\mid & \mid \\
		\bm{w}_1 & \bm{w}_2 \\
		\mid & \mid
	\end{bmatrix}
	= W \begin{bmatrix}
		\mid & \mid \\
		\bm{e}_1 & \bm{e}_2 \\
		\mid & \mid
	\end{bmatrix},\ \quad \text{and} \quad
	\begin{bmatrix}
		\mid \\
		\bm{v}_1 \\
		\mid 
	\end{bmatrix}
	= V \begin{bmatrix}
		\mid \\
		\bm{e}_1 \\
		\mid 
	\end{bmatrix},\]\\
	\item\label{cond:2} $W^{\!T}V = I_{N}$ (biorthogonality),
	\item \label{cond:3} $W^{\!T} Z V=H$.
\end{enumerate}
\end{problem}
Solving this IEP is equivalent to computing a sequence of biorthogonal type I functions and type II MOPs, orthogonal with respect to the discrete inner product \cref{INNER}. We obtain the vectors $\bm{w}_1, \bm{w}_2, \bm{v}_1$ by using the following proposition, which is a direct application of \cref{theorem:VanAsscheThm2}. 
\begin{proposition}\label{prop:initialvectors}
    Assume we have an AT-system (\cref{def:AT-system}) with the positive discrete measures $\mu_{1}$ and $\mu_{2}$ as in \cref{MEASURE}.
    Then, the starting vectors $\bm{w}_1,\bm{w}_2$ and $\bm{v}_1$ for \cref{Prob:IEP} are given as 
    \begin{align*}
        \bm{w}_{1} &= \frac{1}{d_1}\bm{\alpha}_1, \quad \quad \quad \quad \quad \quad  \quad \bm{v}_1 = \begin{bmatrix}
            1 & 1 & \ldots &1
        \end{bmatrix}^T,\\
        \bm{w}_2 &= \frac{1}{d_3}\left(\bm{\alpha}_2 -d_2 \bm{w}_1 \right),
    \end{align*} where
    \begin{equation*}
        d_1 = \sum_{i=1}^{N}\alpha_{1,i},\quad d_2 = \sum_{i=1}^{N}\alpha_{2,i},\quad d_3=\sum_{i=1}^N\left(z_i-\frac{\sum_{j=1}^N z_j\alpha_{1,j}}{\sum_{j=1}^N\alpha_{1,j}}\right) \alpha_{2,i}.
    \end{equation*}
\end{proposition}
\begin{proof}
    It follows from \cref{Prob:IEP} that $\bm{v}_1 = \begin{bmatrix} P_0(z_0) &P_0(z_1) &\cdots &P_0(z_N) \end{bmatrix}^{T}$. From the definition of type II MOPs, $P_0$ is a monic polynomial of degree 0 and thus $P_0 = 1$. Applying \cref{theorem:VanAsscheThm2}, we directly obtain the expressions for $\bm{w}_1$ and $\bm{w}_2$. Likewise, the formulas for $d_1$ and $d_2$ follow immediately. To compute $d_3$, one needs the monic type II polynomial $P_1(x) = x + c$. The constant $c$ can be determined by imposing the orthogonality of $P_1$ with respect to the measure $\mu_{1}$
    \begin{equation*}
        \int_{\R} P_1(x) d\mu_1(x) = \sum_{j=1}^N (z_j+c)\alpha_{1,j}  = 0,
    \end{equation*}
    which gives
    \begin{equation*}
        c = \frac{-\sum_{j=1}^N z_j \alpha_{1,j}}{\sum_{j=1}^N\alpha_{1,j}}.
    \end{equation*}
    Now, $d_3$ can be computed as
    \begin{align*}  
        d_3 = \int_{\mathbb{R}} P_1(x) \, d\mu_2(x) = \sum_{i=1}^{N} (z_i+c) \alpha_{2,i}.
    \end{align*}    
\end{proof}
As already mentioned, the focus is on $r=2$, although everything can be extended to a general $r \geq 2$. For \rev{a} general $r$, the banded recurrence matrix $H$ will have $r$ superdiagonals. Coussement and Van Assche \cite[Theorem $3.2$]{CouVanass05} proved the extension of \cref{theorem:VanAsscheThm2} for \rev{a} general $r$ which can be used to find the initial vectors $\w_1,\w_2,\ldots,\w_r$ to solve the corresponding inverse eigenvalue problem.

\section{Connection with (block) Krylov subspaces}\label{sec:Krylov}
    In this section, we discuss the connection between Krylov subspaces and MOPs. By exploiting this relationship, we specialize the Lanczos-like procedure with multiple starting vectors, as developed by Aliaga et al.~\cite{AliBolFreHer00}, to solve \cref{Prob:IEP}. Given a matrix $M\in\R^{N\times N}$ and a starting vector $\bm v_1\in \R^{N\times N}$, 
     we define the Krylov subspace 
    \begin{equation}
    \label{eq:Krylov-subspace}
        \kryl_m(M,\v_1) = \operatorname{span}\{ \v_1, M\v_1, \dots, M^{m-1} \v_1 \}.
    \end{equation}
    %
    Krylov methods are based on the projection of the matrix $M$ onto $\kryl_m(M,\v_1)$. Since
    the spanning vectors in \eqref{eq:Krylov-subspace} are usually close to being linearly dependent,
    one typically computes a different basis.
    An orthonormal basis $V_m=[\v_1,\v_2,\dots,\v_m]$ of $\kryl_m(M,\v_1)$ (here we assume that $\|\v_1\|=1$) can be computed sequentially by the Arnoldi algorithm, which reduces to the three-term recurrence of Lanczos if $M$ is symmetric. The basis is also nested, meaning that $\{\v_1,\dots,\v_i\}$ is a basis of $\kryl_i(M,\v_1)$ for all $i\leq m$.
    Other approaches include biorthogonal methods, for which the projection onto $\kryl_m(M,\v_1)$ is orthogonal to another Krylov subspace $\kryl_m(M^T,\w_1)$ where $w_1^Tv_1=1$. For instance, the biorthogonal Lanczos algorithm (see, e.g., \cite[Algorithm~7.1]{Saad03}) generates matrices $V_m, W_m \in \mathbb{R}^{N\times m}$, whose columns form bases of $\kryl_m(M,\v_1)$ and $\kryl_m(M^T,\w_1)$, respectively, and satisfy the biorthogonality condition $W_m^T V_m = I_m$. The columns of $V_m$ and $W_m$ satisfy three-term recurrence relations, and the projection $W^T_mMV_m$ is tridiagonal.

    When considering multiple starting vectors, one can use block variants of Krylov subspaces. Given two starting vectors $\w_1,\w_2\in\R^N$ (for the general case see, e.g., \cite{AliBolFreHer00,FroLuSzy17}), we define the block Krylov subspace as
    \begin{equation}
    \label{eq:block-Krylov-subspaces}
        \kryl_m^{\square}(M,[\w_1,\w_2]) 
        =
        \begin{cases}
            \kryl_k(M,\w_1) + \kryl_k(M,\w_2) 
            \quad\quad\text{if $m=2k$},\\
            \kryl_{k+1}(M,\w_1) + \kryl_k(M,\w_2) 
            \quad\text{if $m=2k+1$}.
        \end{cases} 
    \end{equation}
    For notational convenience, the elements of $\kryl_m^{\square}(M,[\w_1,\w_2])$ are vectors. In other common notations, 
    the block Krylov subspace consists of block vectors; see, e.g.,~\cite{Lund18}. However, the latter notation makes it harder to deal with odd indices, as is needed \rev{in our case.} 
    Classical Krylov methods, such as the Arnoldi and Lanczos algorithms, can be extended to the block setting~\cite{FroLuSzy17}. 
    Aliaga et al.~\cite{AliBolFreHer00} study a biorthogonal Lanczos-type algorithm for computing a pair of biorthogonal bases of two block Krylov subspaces with a possibly different number of starting vectors, and the recurrence matrix is banded, with the lower and upper bandwidths depending on the number of starting vectors of each subspace. \cref{alg:biorthogonal-Lanczos-IEP} gives a specialized version for solving \cref{Prob:IEP}.

    Krylov subspaces are profoundly linked with polynomials. For the standard case, every element $\v\in\kryl_m(M,\v)$ can be expressed as $\v = p_{m-1}(M)\v_1$, where $p_{m-1}$ is a polynomial of degree at most $m-1$. Moreover, the orthonormal basis generated by the Arnoldi algorithm corresponds to a system of orthogonal polynomials associated with a suitable inner product; see, e.g., \cite{GolMeu10}. On the other hand, for $\w\in\kryl_m^\square(M,[\w_1,\w_2])$, it follows from \cref{eq:block-Krylov-subspaces} that $\w = p(A)\w_1 + q(A) \w_2$, where $p$ and $q$ are polynomials of appropriate degree. 

    Alqahtani and Reichel~\cite{AlqRei18} show that, for a symmetric $M$, the Krylov subspaces $\kryl_m(M,\v_1)$, $\kryl_m^\square(M^T,[\w_1,\w_2])$ are linked to an unweighted, indefinite inner product and the associated MOPs of type II and I, respectively. Here, we do the opposite: given a family of inner products and the associated MOPs, we link them with the appropriate Krylov subspaces and biorthogonal bases.

    The link between inverse eigenvalue problems and Krylov subspace methods is not new in the literature. For example, the Lanczos or Arnoldi algorithms can be used to solve an IEP with orthonormal eigenvectors and tridiagonal or Hessenberg structure~\cite{BolGol87}, while rational Krylov methods can be used for IEPs with a pencil formulation~\cite{NielVanbVand19,VanbVanbVand22}.
    The following result shows the connection between the IEP in \cref{Prob:IEP} and Krylov subspaces.
    \begin{proposition}
    \label{prop:link-IEP-Krylov}
        Let $Z=\diag(z_1,\dots,z_N)$, $\w_1,\w_2$ and $\v_1$ be as in \cref{Prob:IEP}. Suppose that the IEP admits a solution, i.e., $V=[\v_1,\dots,\v_N]$, $W=[\w_1,\dots,\w_N]$, and $H$ with the structure~\cref{eq:Hessenberg-recurrence-matrix} such that $W^T Z V = H$ and $W^T V = I_N$. Furthermore, assume that $d_n \neq 0$ (cf.~\cref{rem:nonzero-second-superdiagonal}) for $n=2,\dots,N-1$. Then, for $m\leq N$, the sets 
        $\{\v_i\}_{i=1}^m$ and $\{\w_i\}_{i=1}^m$
        form a basis of $\kryl_\dimk(Z,\v_1)$ and $\kryl_m^\square(Z,[\w_1,\w_2])$, respectively.
    \end{proposition}
    \begin{proof}
        Since $W^T V =I_N$, the matrices $V$ and $W$ are both invertible. Hence, both the sets $\{\v_1,\dots,\v_m\}$ and $\{\w_1,\dots,\w_m\}$ form a system of linearly independent vectors. To prove both statements, it suffices to show that $\v_n\in\kryl_n(Z,\v_1)$ and $\w_n\in\kryl_n^\square (Z,[\w_1,\w_2])$ for all $n=1,\dots,N$, since then $\{ \v_1,\dots,\v_m \}\subset\kryl_m(Z,\v_1)$ and $\{ \w_1,\dots,\w_m \}\subset\kryl_m^\square(Z,[\w_1,\w_2])$, and we conclude using the fact that the dimension of both subspaces is at most $m$.
        \begin{enumerate}
            \item Let us show that $\v_n\in\kryl_n(Z,\v_1)$ by induction. By definition, we have that $\v_1\in \kryl_1(Z,\v_1)$. From the relation $ZV=VH$ and the structure of $H$ in~\cref{eq:Hessenberg-recurrence-matrix}, we get that
            \begin{equation}
            \label{eq:vector-recurrence-v}
                \v_{n+1} = Z\v_n -b_n\v_n -c_n\v_{n-1} - d_n\v_{n-2}\quad \text{for $n=1,\dots,N-1$},
            \end{equation}
            where we assume $\v_{-1},\v_0=0$ and $c_0,d_0,d_1=0$. Hence, since $\v_n,\v_{n-1},\v_{n-2}\in\kryl_n(Z,\v_1)$ by inductive hypothesis, and $Z\v_n\in Z\kryl_n(Z,\v_1)\subset\kryl_{n+1}(Z,\v_1)$, we have that
            \begin{equation*}
                \v_{n+1}\in\operatorname{span}\{Z\v_n,\v_n,\v_{n-1},\v_{n-2}\} \subset \kryl_{n+1}(Z,\v_1),
            \end{equation*}
            for all $n=1,\dots,N-1$.
            \item Now we show that $\w_n\in\kryl_n^\square (Z,[\w_1,\w_2])$. By definition, we have that $\w_1\in \kryl_1^\square (Z,[\w_1,\w_2])$ and $\w_2\in\kryl_2^\square (Z,[\w_1,\w_2])$. In this case, we use the relation $ZW = WH^T$ and get that for $n=1,\dots,N-1$,
            \begin{equation}
            \label{eq:vector-recurrence-w}
                \w_{n+1} = \frac{1}{d_n}\left(Z\w_{n-1}  - c_{n-1}\w_n - b_{n-2}\w_{n-1} - \w_{n-2}\right).
            \end{equation}
            By inductive hypothesis, we have that $\w_n,\w_{n-1},\w_{n-2}\in\kryl_n^\square(Z,[\w_1,\w_2])$. Therefore, from the definition of the block Krylov subspace, we get that $Z \w_{n-1}\in Z\kryl_{n-1}^\square(Z,[\w_1,\w_2])\subset \kryl_{n+1}^\square(Z,[\w_1,\w_2])$.
            Thus,
            \begin{equation*}
                \w_{n+1}\in\operatorname{span}\{ Z\w_{n-1}, \w_n,\w_{n-1},\w_{n-2} \}\subset \kryl_{n+1}^\square(Z,[\w_1,\w_2]).
            \end{equation*}
            
        \end{enumerate}
    \end{proof}
\subsection{A biorthogonal Lanczos algorithm}

    By \cref{prop:link-IEP-Krylov}, a solution of \cref{Prob:IEP} consists of a recurrence matrix $H$ and a pair $V$ and $W$ of biorthogonal bases which are nested for the corresponding Krylov subspaces. We consider \cref{alg:biorthogonal-Lanczos-IEP} for solving \cref{Prob:IEP}.

    \begin{algorithm}[ht!]
        \caption{\texttt{IEP\_KRYL}}
        \label{alg:biorthogonal-Lanczos-IEP}
        \begin{algorithmic}[1]
            \REQUIRE{Diagonal matrix $Z=\diag(z_1,\dots,z_N)$, vectors $\bm w_1,\bm w_2, \bm v_1\in\R^n$ such that $\bm w_1^T \bm v_1=1$, $\bm w_2^T \bm v_1 = 0$, $\bm w_2^T Z \bm v_1 = 1$.}
            \ENSURE{$W=[\bm w_1,\dots,\bm w_N]$, $V=[\bm v_1,\dots,\bm v_N]\in \R^{N\times N}$ such that $W^T V = I_N$, and entries $b_0,\dots,b_{N-1}$, $c_1,\dots,c_{N-1}$, $d_2,\dots,d_{N-1}$ of the banded Hessenberg matrix $H$ as in \cref{eq:Hessenberg-recurrence-matrix} such that $W^T Z V = H$\rev{.}}
            \STATE{$b_0 = \bm w_1^T Z \bm v_1$} \label{line:alg1-first}
            \STATE{$\bm v_2 = Z \bm v_1 - b_0\v_1$} \label{line:v2}
            \STATE{$b_1 = \w_2^T Z \v_2$, $c_1 = \w_1^T Z \v_2$}
            \STATE{$\v_3 = Z \v_2 - b_1\v_2 - c_1 \v_1$}
            \STATE{$\hat{\w}_3 = Z\w_1 - c_1\w_2-b_0\w_1$}
            \STATE{$d_2 = \w_1^T Z \v_3$, $\w_3 = \hat{\w}_3 / d_2$}\label{line:alg1-last-before-loop}
            \FOR{$n = 3 : N-1$}
                \STATE{$b_{n-1} = \w_n^T Z \v_n$, $c_{n-1} = \w_{n-1}^T Z \v_n$} \label{line:alg1-first-loop}
                \STATE{$\v_{n+1} = Z \v_n-b_{n-1}\v_n - c_{n-1}\v_{n-1} - d_{n-1}\v_{n-2}$}\label{line:alg1-short-recurrence-v}
                \STATE{$\hat{\w}_{n+1} = Z \w_{n-1} - c_{n-1}\w_n - b_{n-2}\w_{n-1} - \w_{n-2}$}\label{line:alg1-short-recurrence-w}
                \STATE{$d_n = \w_{n-1}^T Z \v_{n+1}$}\label{line:breakdown}
                \STATE{$\w_{n+1} = \hat{\w}_{n+1} / d_n$}\label{line:alg1-last-loop}
            \ENDFOR
            \STATE{$b_{N-1} = \w_N^T Z\v_N$, $c_{N-1} = \w_{N-1}^T Z \v_N$}
        \end{algorithmic}
    \end{algorithm}    
     
    By construction, the outputs satisfy the recurrence relations $ZV=VH$, $ZW=WH^T$, along with the biorthogonality conditions $\w_{n+1}^T \v_k=\w_k^T \v_{n+1} = 0$, for each $n$ and $k=n-2,n-1,n$. To show this for all $n$ and $k$, thus obtaining $W^TV=I_N$, one can use an induction argument, as done in \cite[Proposition~7.1]{Saad03} for the classical biorthogonal Lanczos algorithm. 
    Moreover, $H$ has the structure described in~\cref{eq:Hessenberg-recurrence-matrix}. Since $\w_{n-2},\w_{n-1},\w_n$ are all orthogonal to $\v_{n+1}$, \cref{line:breakdown} can be replaced by the assignment $d_n = {\w}_{n+1}^T \v_{n+1}$, which is more common in the literature~\cite{Saad03}.
    The algorithm consists of a few vector sums, Euclidean inner products and scalar multiplications per iteration, for which the cost is $O(N)$. Thus, the overall complexity is $\mathcal{O}(N^2)$.
    
    If the algorithm terminates successfully, we can guarantee uniqueness.
    \begin{proposition}
    \label{prop:uniqueness-IEP}
        If no breakdown occurs in \cref{alg:biorthogonal-Lanczos-IEP}, then \cref{Prob:IEP} admits a unique solution. 
    \end{proposition}
    \begin{proof}
        Let $V$, $W$, $H$ be any solution, and denote by $\v_1,\dots,\v_N$, $\w_1,\dots,\w_N$ the columns of $V$ and $W$, respectively, and by $b_0,\dots,b_{N-1}$, $c_1,\dots, c_{N-1}$, $d_2,\dots, d_{N-1}$ the entries of $H$ as in \cref{eq:Hessenberg-recurrence-matrix}. Since $H = W^T Z V$ and $H_{i,j} = \w_i^T Z \v_j$, we get that
    \begin{equation}
    \label{eq:coefficients-formula}
        b_{n-1} = \w_n^T Z \v_n,\quad c_{n-1} = \w_{n-1}^T Z \v_n,\quad d_n = \w_{n-1}^T Z \v_{n+1}.
    \end{equation}
        Moreover, because of the identities $ZV=VH$ and $ZW=WH^T$, the basis vectors in the columns of $V$ and $W$ satisfy the recurrence relations \cref{eq:vector-recurrence-v} and \cref{eq:vector-recurrence-w}, respectively, assuming that $\v_{-1},\v_0,\w_{-1},\w_0=0$ and $c_0,d_0,d_1=0$. Hence, by induction, it follows that the solution coincides with the outcome of \cref{alg:biorthogonal-Lanczos-IEP}. 
    \end{proof}

    \subsection{Numerical aspects}
    \label{subsec:numerical-aspects}

    Biorthogonal Krylov methods often suffer from numerical issues. 
    Here, we discuss how they affect the solution of \cref{Prob:IEP} via \cref{alg:biorthogonal-Lanczos-IEP} and possible strategies to overcome them.

    \rev{\textbf{Look-Ahead.}} A breakdown occurs in \cref{line:breakdown} if $d_n$ equals zero. \Cref{rem:nonzero-second-superdiagonal} ensures that this cannot happen when the IEP is associated with some classes of well-defined MOPs. However, this does not exclude the possibility of a near-breakdown, which arises when $d_n$ is small. A popular workaround for this issue is the look-ahead technique which is described, for instance, by Freund~\cite{Fre94} for the standard biorthogonal Lanczos method and extended by Aliaga et al.~\cite{AliBolFreHer00} to the case of multiple starting vectors. However, this approach changes the sparsity pattern of the recurrence matrix by introducing nonzero entries outside the prescribed bandwidth, thereby perturbing the required structure in \cref{Prob:IEP}. We also found that the near-breakdown, unlike other phenomena, is not the main source of inaccuracy in our tests as we will explain in section \ref{sec:numericalexperiments}. Although we implemented look-ahead techniques in our algorithms, they did not lead to improved accuracy. Therefore, we omit them from both the discussion and the final algorithmic formulations.

\rev{\textbf{Orthogonalization Procedures.}}  In order to get $\v_{n+1}$ and $\w_{n+1}$, we only need to orthogonalize $Z\v_n$ and $Z\w_{n-1}$ against $\w_{n-2}, \w_{n-1}, \w_n$ and $\v_{n-2},\v_{n-1},\v_n$, respectively. In fact, in exact arithmetic, $Z\v_n$ is orthogonal to $\kryl_{n-3}^\square (Z,[\w_1,\w_2])$ and $Z\w_{n-1}$ is orthogonal to $\kryl_{n-3}(Z,\v_1)$. However, numerically, the quantities $\w_i^T \v_{n+1}$ and $\w_{n+1}^T\v_j$ may not be small. This phenomenon affects most Krylov subspace methods and is referred to as loss of orthogonality. To address this issue, one could perform a full biorthogonalization instead of relying on the short recurrence relation. 
    To further ensure the biorthogonality, we can also apply reorthogonalization to the computed bases $W$ and $V$, at the cost of another long-term recurrence formula for each step. For further details and strategies like selective reorthogonalization, see~\cite{MorNic11}.

\textbf{Normalization Strategies.}
    A feature of \cref{alg:biorthogonal-Lanczos-IEP} is that it returns $H$ with lower subdiagonal entries equal to $1$. This is because $\v_{n+1}$ is not scaled after its orthogonalization at \cref{line:alg1-short-recurrence-v}, while the scaling is concentrated on $\w_{n+1}$ at \cref{line:alg1-last-loop}, as required in \cref{Prob:IEP}. These ones on the subdiagonal correspond to the choice of monic type II MOPs as discussed in \cref{sec:steplinerecurrence}. Numerically, the resulting bases are typically severely ill-conditioned, implying that neither the biorthogonality condition $W^T V = I_N$ nor $W^T Z V = H$ will be satisfied \rev{numerically}.
    Different scaling strategies, discussed further on, can be used to mitigate this high conditioning number. 
    For example, by imposing that all basis vectors have unit norm, the biorthogonality relation reduces to $W^TV=\Sigma$, where $\Sigma$ is a diagonal matrix, though not necessarily the identity.\footnote{Throughout this paper, we say that a pair of biorthogonal bases satisfies $W^TV=I_N$, hence excluding the case where $W^TV$ is just a diagonal matrix.} As a consequence, $V$ and $W$ will satisfy the two recurrence relations
    \begin{equation}
    \label{eq:double-recurrence}
        ZV = VH^V, \quad ZW = W H^W,
    \end{equation}
    where $(H^W)^T$ and $H^V$ have the same sparsity pattern as in \cref{eq:Hessenberg-recurrence-matrix} but without the ones on the subdiagonal. This normalization is used in \cref{alg:biorthogonal-Lanczos-IEP-generalized}. 
    \rev{
    Note that, after a full orthogonalization, the upper triangular part of $H^V$ and $H^W$ can contain small nonzero elements. However, if the bases are accurately computed, the entries of $H^V$ and $(H^W)^T$ outside the sparsity pattern in \cref{eq:Hessenberg-recurrence-matrix} are negligible and can be set to $0$.
    }
  
\rev{Let us discuss the link between the matrices in \cref{eq:double-recurrence} 
and the output of \cref{alg:biorthogonal-Lanczos-IEP}.} We can get \rev{biorthogonal} bases through a diagonal scaling. More specifically, if $W^TV=\Sigma_1 \Sigma_2$, for any diagonal matrices $\Sigma_1$ and $\Sigma_2$, then $\widetilde{W} = W \Sigma_1^{-1}$ and $\widetilde{V} = V \Sigma_2^{-1}$ satisfy $\widetilde{W}^T \widetilde{V} = I_N$. Furthermore, we get that $\Sigma_\rev{2} H^V \Sigma_{\rev{2}}^{-1} = \Sigma_{\rev{1}}^{-1} (H^W)^T \Sigma_{\rev{1}} = H$ satisfies $\widetilde{W}^T Z \widetilde{V} = H$. \rev{In particular, if $\Sigma_1=I_N$, then $\widetilde{W}=W$ and $H=(H^W)^T$.} 
    
    \rev{The ones on the subdiagonal of $H$ can be obtained via an additional diagonal scaling. Let $a_i=H_{i+1,i}$ for $i=1,\dots,N-1$ be the subdiagonal entries of $H$, define
    \begin{equation}
        \label{eq:diagonal-scaling-ones}
        D=\operatorname{diag}(1,a_1,a_1a_2,\dots,a_1a_2\cdots a_{N-1}).
    \end{equation}
    Then $\widehat H = D^{-1}HD$ is such that $\widehat{H}_{i+1,i}=1$ for all $i$.
    The scaling can be applied to the bases, $\widehat{W} = \widetilde{W}D$ and $\widehat{V}= \widetilde VD^{-1}$, so that $H,W,V$ solve \cref{Prob:IEP}. 
   One should be careful, however, since the factors in \cref{eq:diagonal-scaling-ones} tend to be extremely large or small if, for example, the values $a_i$ have similar magnitude, as a result $D$ can become very ill-conditioned. This makes the scaling required for $\widehat{W}$ and $\widehat{V}$ infeasible in practice, which makes sense, as we revert back to the original ill-conditioned bases. The recurrence matrix is less affected by this issue, since the scaling acts only on the nonzero entries of the sparsity pattern \cref{eq:Hessenberg-recurrence-matrix}, i.e., $\widehat H_{i+1,i} = 1,$ $\widehat H_{i,i} = H_{i,i}$, and
    \begin{equation}
    \label{eq:recurrence-scaling-monic-entries}
        \widehat H_{i,i+1} = a_i H_{i,i+1},\quad \widehat H_{i,i+2} = a_i a_{i+1} H_{i,i+2}.
    \end{equation}
    In this case, the factors $a_i$ and $a_ia_{i+1}$ for the first and second superdiagonal do not typically have a prohibitive difference in magnitude. }

\subsection{An algorithm with normalized vectors and (re-)orthogonalization}
    In \cref{alg:biorthogonal-Lanczos-IEP-generalized} we consider a process including the aforementioned techniques\rev{, which is normalizing the basis vectors, using partial or full (re-)orthogonalization, and a final scaling to get biorthogonal bases and the ones on the subdiagonal of the recurrence matrix, corresponding to monic type II MOPs}. For a matrix $M$, we use the notation $M_{i_1:i_2,j}=[M_{i_1,j}, M_{i_1+1,j},\dots,M_{i_2-1,j},M_{i_2,j}]^T$.

    \begin{algorithm}[h!]
        \caption{\texttt{IEP\_KRYLREORTH}\label{alg:biorthogonal-Lanczos-IEP-generalized}}
        \begin{algorithmic}[1]
            \REQUIRE{Diagonal matrix $Z=\diag(z_1,\dots,z_N)$, vectors $\bm w_1,\bm w_2, \bm v_1\in\R^n$ such that $\bm w_1^T \bm v_1=1$, $\bm w_2^T \bm v_1 = 0$, $\bm w_2^T Z \bm v_1 \neq 1$, and $\|\v_1\| = \|\w_1\| = \|\w_2\| = 1$, orthogonalization strategy \texttt{partial} or \texttt{full}\rev{.}}
            \ENSURE{$W=[\bm w_1,\dots,\bm w_N]$, $V=[\bm v_1,\dots,\bm v_N]$, and $\Sigma =  \diag(\sigma_1,\dots,\sigma_N)$ such that $W^TV=\Sigma$ and $\|\v_n\|=\|\w_n\|=1$ for all $n$, and recurrence matrices $H^V, H^W$  such that $Z V = V H^V$ and $Z W = W H^W$.}
            \STATE{$ \hat{\v}_2 = Z \v_1 - H^V_{1,1}\v_1$, \quad $H^V_{2,1} = \|\hat{\v}_2\|$,\quad $\v_2 = \hat{\v}_2 / H^V_{2,1}$,\quad $\sigma_2 = \w_2^T \v_2$}
            \FOR{$n = 2 : N-1$}\label{line:alg2-internal-loop-start}
                \STATE{$\hat{\v}_{n+1} = Z \v_n$,\quad $\hat{\w}_{n+1} = Z \w_{n-1}$}
                \STATE{Orthogonalization: set $r=\max(1,n-2)$ for \texttt{partial}, $r=1$ for \texttt{full}\label{line:alg2-orthogonalization-strategy}}
                \STATE{Set $V_n = [\v_r,\dots,\v_k]$, $W_n = [\w_r,\dots,\w_n]$, and $\Sigma_n = \diag(\sigma_r,\dots,\sigma_n)$ \label{line:alg2-auxiliary}}
                \STATE{$H^V_{r{:}n,n} = \Sigma_n^{-1} W_n^T \hat{\v}_{n+1} $,\quad $H^W_{r{:}n,n-1} = \Sigma_n^{-1} V_n^T\hat{\w}_{n+1}$ \label{line:alg2-recurrence-coefficients}}
                \STATE{$\hat{\v}_{n+1} = \hat{\v}_{n+1} - V_n H^V_{r:n,n} $,\quad $\hat{\w}_{n+1} = \hat{\w}_{n+1} - W_n H^W_{r:n,n-1} $}
                \STATE{$H^V_{n+1,n} = \|\hat{\v}_{n+1}\|$,\quad $H^W_{n+1,n-1} = \|\hat{\w}_{n+1}\|$}
                \STATE{$\v_{n+1} = \hat{\v}_{n+1} / H^V_{n+1,n}$,\quad $\w_{n+1} = \hat{\w}_{n+1} / H^W_{n+1,n-1}$}
                \STATE{For reorthogonalization, \textbf{go to} line \ref{line:alg2-internal-loop-start} and update $H^V$ and $H^W$}
                \STATE{$\sigma_{n+1} = \w_{n+1}^T \v_{n+1}$}
            \ENDFOR
            \STATE{Repeat lines \ref{line:alg2-orthogonalization-strategy}, \ref{line:alg2-auxiliary},  \ref{line:alg2-recurrence-coefficients} with $n=N$ to get $H^V_{1{:}N,N}$, $H^W_{1:N,N-1}$, $H^W_{1:N,N}$}
            \STATE{\rev{Truncation (for \texttt{full}): $H^V_{i,j} = 0$ for $j>i+2$, $H^W_{i,j} = 0$ for $j>i+1$}}
            \STATE{\rev{Biorthogonality: set $\widetilde{V}=V\Sigma^{-1}$, $\widetilde{W}=W$, and $H=(H^W)^T$}}
            \STATE{\rev{Monic type II: set $D$ as in \eqref{eq:diagonal-scaling-ones}, then $\widehat{H}=D^{-1}HD$, $\widehat{W}=\widetilde{W}D$, $\widehat{V}=\widetilde{V}D^{-1}$}}
        \end{algorithmic}
\end{algorithm}

If partial orthogonalization is used, the complexity  is still $\mathcal{O}(N^2)$ as for \cref{alg:biorthogonal-Lanczos-IEP}, both with and without reorthogonalization.
If we use a full orthogonalization, the algorithm becomes much more expensive, yielding a cost of $\mathcal{O}(n\,N)$ at each step $n$, and therefore an overall complexity of $\mathcal{O}(N^3)$. However, the increased accuracy justifies the higher cost, as illustrated by the experiments in \cref{sec:numericalexperiments}. \rev{Moreover, in many practical applications, the number of required polynomials is limited; see, e.g.,~\cite{LauMastVanAVanDooren24}, implying a modest increase in computational cost.}

\subsection{Moment matrix formulation}
To solve the IEP \cref{Prob:IEP}, we look for nested biorthogonal subspaces
for given vector spaces.
Besides the block Krylov method as described before, these subspaces can be related to moment matrices.

The moment matrix $ {M}_{N} \in \mathbb{R}^{N \times N} $, associated with two sets of vectors 
$ X = \{\bm{x}_i\}_{i=1}^{N} $ and $ Y = \{\bm{y}_i\}_{i=1}^{N} $ in vector spaces $ \mathcal{X} $ and $ \mathcal{Y} $, respectively, is defined as  
\begin{equation}\label{MO}
{M}_{N} = \left( \bm{y}_i^T \bm{x}_j \right)_{i,j=1}^N.
\end{equation}
The following lemma states how to orthogonalize $ X $ and $ Y $ by factoring $ {M}_{N} $.
\begin{lemma}[Biorthonormal vectors via moment matrix factorization \cite{Dav63}]
Let $ X = [\bm{x}_1 \; \bm{x}_2 \; \dots \; \bm{x}_N] $ and $ Y = [\bm{y}_1 \; \bm{y}_2 \; \dots \; \bm{y}_N] $ be full-rank matrices in some vector spaces $ \mathcal{X} $ and $ \mathcal{Y} $, respectively. Let $ {M}_{N} $  be the associated moment matrix \cref{MO}. The (non-pivoted) LR factorization of $ {M}_{N} $ is given by ${M}_{N}  = LR$, where $ L $ is a lower triangular matrix and $ R $ is an upper triangular matrix. Define 
\[
V = [\bm{v}_1 \; \bm{v}_2 \; \dots \; \bm{v}_N] = X R^{-1}, \quad W = [\bm{w}_1 \; \bm{w}_2 \; \dots \; \bm{w}_N] = Y L^{-\!T}.
\]  
Then, for $i = 1, \dots, N$, we have
\[
\operatorname{span} \{ \bm{v}_1, \dots, \bm{v}_i \} = \operatorname{span} \{ \bm{x}_1, \dots, \bm{x}_i \},  
\]  
\[
\operatorname{span} \{ \bm{w}_1, \dots, \bm{w}_i \} = \operatorname{span} \{ \bm{y}_1, \dots, \bm{y}_i \}, 
\]  
and the vectors are biorthonormal  
\[
W^{\!T}V= L^{-1} {M}_{N} R^{-1} = I_N.
\]  
\end{lemma}
Consider $ Z = \diag(\{z_i\}_{i=1}^{N}) $ with distinct diagonal entries, and
the matrices
\begin{equation*}
    K_V = \begin{bmatrix} \bm{v}_{1} & Z\bm{v}_{1} & \dots & Z^{N-1}\bm{v}_{1} \end{bmatrix},
\end{equation*}
 \begin{equation*}
    K_W 
        =
        \begin{cases}
        \vspace{.2cm}
           \begin{bmatrix} \bm{w}_{1} ~\bm{w}_{2} & Z\bm{w}_{1}~ Z\bm{w}_{2}& \dots & Z^{k-1}\bm{w}_{1}~  Z^{k-1}\bm{w}_{2}\end{bmatrix},
            \quad\quad\quad \quad~\text{if $N=2k$},\\
            
            \begin{bmatrix} \bm{w}_{1} ~\bm{w}_{2} & Z\bm{w}_{1}~ Z\bm{w}_{2}& \dots & Z^{k-1}\bm{w}_{1}~  Z^{k-1}\bm{w}_{2} &Z^{k}\bm{w}_{1}\end{bmatrix},
            \quad \text{if $N=2k+1$}.
        \end{cases} 
    \end{equation*}
    The columns of $K_V$ and $K_W$ form nested bases of the Krylov subspcaes $\kryl_m(Z,\v_1)$ and $\kryl_m^\square(Z,[\w_1,\w_2])$, respectively, for all $m=1,\dots,N$. Thus, we refer to the columns of $K_V$, $K_W$ as Krylov bases. The moment matrix can now be constructed as ${M}_{N}=K_W^{\!T} K_V$, and it represents the moment matrix for a sequence of biorthonormal MOPs. A direct way to compute the recurrence coefficients of a set of biorthonormal MOPs is to evaluate the LR factorization of $\mathcal{M}_{N}$. In finite precision, two issues are encountered with this procedure:
\begin{itemize}
    \item The LR factorization without pivoting is numerically unstable,
    \item The matrix ${M}_{N}$ tends to be very ill-conditioned.
\end{itemize}
The use of modified moments~\cite{Gautschi04} could improve the conditioning of the moment matrix; however, we have not explored this approach in this setting.

\section{Core transformation algorithm}\label{Sec5}
\noindent
This section presents an algorithm based on core transformations to solve the inverse eigenvalue problem described in \cref{Prob:IEP}. The algorithm obtains the unique recurrence matrix $H$ by applying a sequence of non-unitary similarity transformations $P^{-1}_k\ldots P_1^{-1}ZP_1\ldots P_k$ on $Z$.\footnote{In the literature, these transformations are typically referred to as core transformations~\cite{CampsMeerbVandeb2019}.} These matrices $P_i$, called eliminators, are upper or lower triangular with a $2\times2$ active part acting on two consecutive rows or columns. More precisely, they are identity matrices with a $2 \times 2$ lower or upper matrix, denoted by \\
\tikzexternaldisable
\vspace{-1em}
\begin{equation}
\begin{aligned}
&\begin{array}{c}
    \tikzmarknode[anchor=base]{a1}{ } \\
    \tikzmarknode[anchor=base]{a2}{ }
\end{array}
= \begin{bmatrix}
    1 & 0 \\
    \tau & 1
\end{bmatrix}\quad \text{and}
\quad
\begin{tikzpicture}[remember picture,overlay]
    \draw[-latex] 
        ([xshift=-0.5ex,yshift=0.4ex]a1.west) 
        -- ++(-0.2,0)  
        |- ([xshift=-0.2ex, yshift=0ex]a2.west);
\end{tikzpicture}
\quad
\quad
\begin{array}{c}
    \tikzmarknode[anchor=base]{b2}{ } \\
    \tikzmarknode[anchor=base]{b1}{ }
\end{array}
= \begin{bmatrix}
    1 & \tau \\
    0 & 1
\end{bmatrix},
\quad
\begin{tikzpicture}[remember picture,overlay]
    \draw[-latex] 
        ([xshift=-0.5ex,yshift=0.4ex]b1.west) 
        -- ++(-0.2,0)  
        |- ([xshift=-0.2ex, yshift=0ex]b2.west);
\end{tikzpicture}
\end{aligned}
\label{eq:cores}
\end{equation}
\vspace{-0.1em}
\tikzexternaldisable
\noindent
embedded along the diagonal.
This notation illustrates that a multiple of one row (or column) is added to the next or previous row (or column). For further reading on how eliminators can be used to solve inverse eigenvalue problems, we refer to~\cite{VB21,CampsMeerbVandeb2019,MachVanbVand14}. \cref{alg:coretransformation} provides an outline of the algorithm for a general dimension $N$. To illustrate the procedure, we consider an example with $N = 5$. We first do $N-1$ steps which consist of three parts: \emph{elimination}, \emph{biorthogonalization} and \emph{chasing}. By eliminating elements in $\bm{w}_1, \bm{w}_2$ and $\bm{v}_1$ using eliminators, condition \ref{cond:1} will be satisfied. Hence, we need to find matrices $W^{-1}$ and $V^{-1}$ such that $W^{-1}\begin{bmatrix}
	\bm{w}_1 & \bm{w}_2 
\end{bmatrix} = \begin{bmatrix}
	\bm{e}_1 & \bm{e}_2
\end{bmatrix}$ and $V^{-1}\bm{v}_1 = \bm{e}_1$. 
By using LU factorizations in the biorthogonalization step, we can guarantee that condition \ref{cond:2} is satisfied. The matrices needed in the elimination and biorthogonalization steps are applied to the diagonal matrix $Z$. To achieve the desired banded structure \ref{cond:3} of $H$, a chasing step is performed to remove undesired elements. The final step applies a diagonal similarity transformation on $H$ to obtain the ones on the subdiagonal.

\begin{algorithm}[h!]     
	\caption{\texttt{IEP\_CORE}}
    \label{alg:coretransformation}
	\begin{algorithmic}[1] 
		\REQUIRE{Diagonal matrix $Z=\diag(z_1,\dots,z_N)$, vectors $\bm w_1,\bm w_2, \bm v_1\in\R^N$ such that $\bm w_1^T \bm v_1=1$, $\bm w_2^T \bm v_1 = 0$ and $\bm w_2^T Z \bm v_1 = 1$.}
		\ENSURE{$W=[\bm w_1,\dots,\bm w_N], V=[\bm v_1,\dots,\bm v_N]\in \R^{N\times N}$ such that $W^T V = I$, and $H$ banded Hessenberg matrix as in \cref{eq:Hessenberg-recurrence-matrix} such that $W^T Z V = H$.} \\
        {\texttt{\% Step 1}}
		\STATE{Compute eliminators $W_{(1,N)}^{-1}, W_{(1,N-1)}^{-1}, W_{(2,N)}^{-1}$ and $V_{(1,N)}^{-1}$}  
		\STATE{Compute $\widehat{W}_1^{-1} = W^{-1}_{(1,N)}W^{-1}_{(1,N-1)}W^{-1}_{(2,N)}$ and $\widehat{V}_1^{-1} = V^{-1}_{(1,N)}$}
		\STATE{$\widehat{W}_{1}^{T}\widehat{V}_{1} = L_1 U_1$} 
		\STATE{$W_{1} = \widehat{W}_{1}L_1^{-1}$ and $V_{1} = \widehat{V}_{1}U_1^{-1}$} 
		\STATE{$W^T = W_{1}, V = V_{1}$ and $H^{(1)} = W^T Z V$}\\
        {\texttt{\% Step 2 to $N-2$}}
		\FOR{$i = 2:N-2$}
		\STATE{Compute eliminators $W_{(1,N-i)}, W_{(2,N-i-1)}$ and $V_{(1,N-i-1)}$}
		\STATE{$\widehat{W}_{i}=W_{(1,N-i)}W_{(2,N-i-1)}$ and $\widehat{V}_{i} = V_{(1,N-i-1)}$}
		\STATE{$\widehat{W}_{i}^{T}\widehat{V}_{i} = L_i U_i$} 
		\STATE{$W_{i} = \widehat{W}_{i} L_i^{-T}$ and $V_{i} = \widehat{V}_{i}U_i^{-1}$} 
		\STATE{$W^T=W^T_{i}W^T$, $V = V V_{i}$, $H^{(i)} = W_i^TH^{(i-1)}V_i$}
		\FOR{$j = 0:i-1$}
		\STATE{Eliminate $H^{(i)}_{N-i+j+1, N-i+j}$}   {\quad \;\;\quad \texttt{\% Chase lower bulge}}
		\STATE{Eliminate $H^{(i)}_{N-i+j-1,N-i+j+2}$}  {\quad \texttt{\% Chase upper bulge}}
		\ENDFOR
		\ENDFOR\\
        {\texttt{\% Step $N-1$}}
        \STATE{Compute ${W}^T_{N-1}$ and $V_{N-1} = {W}^{-1}_{N-1}$} 
        \STATE $W^T=W^T_{N-1}W^T$, $V = V V_{N-1}$, $H^{(N-1)} = W_{N-1}^TH^{(N-2)}V_{N-1}$
		\FOR{j=0:(N-3)} 
            \STATE{Eliminate $H^{(N-1)}_{j+3,j+1}$} {\qquad \qquad \qquad \qquad \texttt{\% Chase lower bulge}}
        \ENDFOR\\
        {\texttt{\% Step $N$}}
        \STATE{Apply $S_1,S_2,S_3$ such that $\bm{e}_1 = S_1 \bm{w}^{(N-1)}_1,\bm{e}_2 = S_2 \bm{w}^{(N-1)}_2$ and $\bm{e}_1 = S_3 \bm{v}^{(N-1)}_1$} 
        \STATE{Apply scaling $W_s$ such that lower diagonal of $H^{(N)}=W_sH^{(N-1)}W_s^{-1}$ are all ones}
    \end{algorithmic}
\end{algorithm}

\begin{figure}[h!]
	\begin{subfigure}[t]{.33\textwidth}
		\centering
		$\setlength\arraycolsep{1pt}
		\begin{array}{cc|ccccc}
			& & {v}_1^1 & {v}_1^2 & {v}_1^3 & \tikzmarknode{v4}{v_1^4} & \tikzmarknode{v5}{v_1^5} \\
			\hline \Tstrut
			w_2^{1}\; \; &  w_1^{1}&  z_1 & & & & \\
			w_2^{2}\; \; &  w_1^{2}&  & z_2 & & & \\
			\tikzmarknode{w3}{w_2^{3}}\; \; &  w_1^{3}&  & & z_3 & & \\
			\tikzmarknode{w4}{w_2^{4}}\; \; &  w_1^{4}&  & & & z_4 & \\
			\tikzmarknode{w5}{w_2^{5}}\; \; &  w_1^{5}&  & &  &  & z_5 \\
		\end{array}
        \begin{tikzpicture}[remember picture,overlay]
			\draw[-latex] 
			([yshift=.25ex]v4.north) 
			-- ++(0,.16)  
			-| (v5.north);  
			\draw[-latex] 
			([xshift=-0.25ex]w4.west) 
			-- ++(-.16,0)  
			|- ([xshift=-.ex]w5.west);  
			\draw[-latex] 
			([xshift=-1.6ex]w3.west) 
			-- ++(-.16,0)  
			|- ([xshift=-1.35ex]w4.west);
			\draw[-latex] 
			([xshift=-2.85ex]w4.west) 
			-- ++(-.16,0)  
			|- ([xshift=-2.7ex]w5.west); 
		\end{tikzpicture}
		$
		\caption{Step 1: Elim. \& Bior.}
		\label{fig:step0}
	\end{subfigure}%
	\begin{subfigure}[t]{.33\textwidth}
		\centering
		$\setlength\arraycolsep{1pt}
		\begin{array}{cc|ccccc}
			& & {v}_1^1 & {v}_1^2 & \tikzmarknode{v3}{v_1^3} & \tikzmarknode{v4}{v_1^4} & v_1^5 \\
			\hline \Tstrut
			w_2^{1}\; \; &  w_1^{1}&  z_1 & & & & \\
			\tikzmarknode{w2}{w_2^{2}}\; \; &  w_1^{2}&  & z_2 & & & \\
			\tikzmarknode{w3}{w_2^{3}}\; \; &  w_1^{3}&  & & z_3 & \times &\times \\
			\tikzmarknode{w4}{\times} \; \; &  &  & & & \times & \times\\
			\; \; &  &  & &  & \times & \times \\
		\end{array}
        \begin{tikzpicture}[remember picture,overlay]
			\draw[-latex] 
			([yshift=.25ex]v3.north) 
			-- ++(0,.16)  
			-| (v4.north);  
			\draw[-latex] 
			([xshift=-.25ex]w2.west) 
			-- ++(-.16,0)  
			|- ([xshift=-0.ex]w3.west);  
			\draw[-latex] 
			([xshift=-1.6ex]w3.west) 
			-- ++(-.16,0)  
			|- ([xshift=-1.6ex]w4.west); 
		\end{tikzpicture}
		$
		\caption{Step 2: Elim. \& Bior.}
		\label{fig:step1}
	\end{subfigure}%
	\begin{subfigure}[t]{.33\textwidth}
		\centering
		$\setlength\arraycolsep{1pt}
		\begin{array}{cc|ccccc}
			& & {v}_1^1 & {v}_1^2 & \tikzmarknode{v3}{{v}_1^3} & \tikzmarknode{v4}{v_1^4} & \tikzmarknode{v5}{ } 	\\
			\hline \Tstrut
			w_2^{1}\; \; &  w_1^{1}&  z_1 & & & & \\
			w_2^{2}\; \; & w_1^2  &  & z_2 & \times & \times & \thicktimes \\
		      \times \; \; &  &  & & \times & \times &\times \\
			  \tikzmarknode{w4}{ } &  &  & & \times& \times &\times \\
			  \tikzmarknode{w5}{ }&  &  & &  \thicktimes & \times& \times \\
		\end{array}
        \begin{tikzpicture}[remember picture,overlay]
            \draw[-latex,red]
			([yshift=1.8ex]v5.north) 
			-- ++(0,.12)  
			-| ([yshift=-.5ex]v4.north);
            \draw[-latex,red]
			([yshift=1.2ex]v4.north) 
			-- ++(0,.12)  
			-| ([yshift=2.6ex]v5.north); 
			\draw[-latex,red] ([xshift=-1.1ex,yshift=0.4ex]w4.west) 
			-- ++(-.16,0) 
			|- ([yshift=.4ex,xshift=-0.9ex]w5.west);
            \draw[-latex,red] ([xshift=-2.45ex,yshift=0.4ex]w5.west) 
			-- ++(-.16,0) 
			|- ([yshift=.4ex,xshift=-2.25ex]w4.west);
		\end{tikzpicture}
		$
		\caption{Step 2: Chasing}
		\label{fig:step2a}
	\end{subfigure}%
    \\[5pt]
	\begin{subfigure}[t]{.33\textwidth}
		\centering
		$\setlength\arraycolsep{1pt}
		\begin{array}{cc|ccccc}
			& & {v}_1^1 & \tikzmarknode{v2}{{v}_1^2} & \tikzmarknode{v3}{\times} &  & \\
			\hline \Tstrut
			\tikzmarknode{w1}{w_2^{1}}\; \; &  w_1^{1}&  z_1 & & & & \\
			\tikzmarknode{w2}{w_2^{2}}\; \; &  w_1^{2}&  & z_2 & \times & \times &  \\
			\tikzmarknode{w3}{\times} &  &  & & \times& \times &\times \\
			\; \; &  &  & & \times& \times &\times \\
			\; \; &  &  & &  & \times& \times \\
		\end{array}
        \begin{tikzpicture}[remember picture,overlay]
			\draw[-latex] 
			([yshift=.25ex]v2.north) 
			-- ++(0,.16)  
			-| ([yshift=0.30ex]v3.north);  
			\draw[-latex] 
			([xshift=-.25ex]w1.west) 
			-- ++(-.16,0)  
			|- (w2.west);  
			\draw[-latex] 
			([xshift=-1.6ex]w2.west) 
			-- ++(-.16,0)  
			|- ([xshift=-1.9ex]w3.west); 
		\end{tikzpicture}
		$
		\caption{Step 3: Elim. \& Bior.}
		\label{fig:step2b}
	\end{subfigure}%
	\begin{subfigure}[t]{.33\textwidth}
		\centering
		$\setlength\arraycolsep{1pt}
		\begin{array}{cc|ccccc}
			& & \times & \times & \tikzmarknode{v3}{ } &\tikzmarknode{v4}{ } & 	\\
			\hline \Tstrut
			\times\; \; &  \times& z_1 & \times & {\times} & \thicktimes & \\
			\times\; \; &   & &\times & \times & \times & \thicktimes \\
			\tikzmarknode{w3}{ }\; \; & &\tikzmarknode{3l}{ }  & \times& \times & \times &\times \\
			  \tikzmarknode{w4}{ }&  &\tikzmarknode{4l}{ }  & \thicktimes & \times& \times &\times \\
		    &  &  & &  & \times& \times \\
		\end{array}
		\begin{tikzpicture}[remember picture,overlay]
            \draw[-latex,red]
			([yshift=2ex]v4.north) 
			-- ++(0,.16)  
			-| ([yshift=1.5ex]v3.north);  
			\draw[-latex,red] ([xshift=-.8ex,yshift=1.4ex]w3.west) 
			-- ++(-.16,0) 
			|- ([yshift=.7ex,xshift=-0.9ex]w4.west);
		\end{tikzpicture}
		$
		\caption{Step 3: Lower chasing}
		\label{fig:step3a}
	\end{subfigure}%
	\begin{subfigure}[t]{.33\textwidth}
		\centering
		$\setlength\arraycolsep{1pt}
		\begin{array}{cc|ccccc}
			& & \times & \times &  & & 	\\
			\hline \Tstrut
			\times\; \; &  \times & z_1 & \times &\times & \thicktimes& \\
			\times \; \; &  &  & \times & \times & \times & \thicktimes \\
			\; \; &  &  & \times& \times & \times &\times \\
			\; \; &  &  & & \times& \times &\times \\
			\; \; &  &  & &  \thicktimes & \times& \times \\
		\end{array}
		$
		
		\caption{1 lower bulge chased}
		\label{fig:step3b}
	\end{subfigure}%
    \\[5pt]
	\begin{subfigure}[t]{.25\textwidth}
		\centering
		$\setlength\arraycolsep{1pt}
		\begin{array}{cc|ccccc}
			& & \tikzmarknode{v1}{\times} & \tikzmarknode{v2}{\times} &  &  & \\
			\hline \Tstrut
			\tikzmarknode{w1}{\times} \; \; & \times &  z_1 & \times &\times & & \\
			\tikzmarknode{w2}{\times}\; \; &  &  & \times & \times & \times &  \\
			\; \; &  &  &\times & \times& \times &\times \\
			\; \; &  &  & & \times& \times &\times \\
			\; \; &  &  & &  & \times& \times \\
		\end{array}
        \begin{tikzpicture}[remember picture,overlay]
			\draw[-latex]
			([yshift=.5ex]v1.north) 
			-- ++(0,.16)  
			-| ([yshift=0.3ex]v2.north);  
			\draw[-latex] 
			([xshift=-.3ex]w2.west) 
			-- ++(-.16,0)  
			|- ([yshift=.2ex,xshift=.2ex]w1.west);  
		\end{tikzpicture}
		$
		\caption{Step 4: Elim.}
		\label{fig:step4a}
	\end{subfigure}%
	\begin{subfigure}[t]{.25\textwidth}
		\centering
		$\setlength\arraycolsep{1pt}
		\begin{array}{cc|ccccc}
			& & \times & &  &  & \\
			\hline \Tstrut
			\; \; &  \times&  \times & \times &\times & & \\
			\times\; \; &  &  \times& \times & \times & \times &  \\
			\; \; &  &  \thicktimes&\times & \times& \times &\times \\
			\; \; &  &  & & \times& \times &\times \\
			\; \; &  &  & &  & \times& \times \\
		\end{array}
		$
		\caption{Final bulge chase}
		\label{fig:step4b}
	\end{subfigure}%
	\begin{subfigure}[t]{.25\textwidth}
		\centering
		$\setlength\arraycolsep{1pt}
		\begin{array}{cc|ccccc}
			& & \times &  &  &  & \\
			\hline \Tstrut
			\; \; &  \times&  \times & \times &\times & & \\
			\times\; \; &  &  \times& \times & \times & \times &  \\
			\; \; &  &  &\times & \times& \times &\times \\
			\; \; &  &  & & \times& \times &\times \\
			\; \; &  &  & &  & \times& \times \\
		\end{array}
		$
		\caption{Entire step 4}
		\label{fig:step4c}
	\end{subfigure}
	\begin{subfigure}[t]{.24\textwidth}
		\centering
		$\setlength\arraycolsep{1pt}
		\begin{array}{cc|ccccc}
			& & 1 &  &  &  & \\
			\hline \Tstrut
			\; \; &  1&  \times & \times &\times & & \\
			1\; \; &  &  1& \times & \times & \times &  \\
			\; \; &  &  &1 & \times& \times &\times \\
			\; \; &  &  & & 1& \times &\times \\
			\; \; &  &  & &  & 1& \times \\
		\end{array}
		$
		\caption{Step 5: Scaling}
		\label{fig:step5}
	\end{subfigure}%
	\caption{Visualization of different steps in \cref{alg:coretransformation} for $N=5$. To simplify the notation, we do not display the transformations for the biorthogonalization.}
\end{figure}
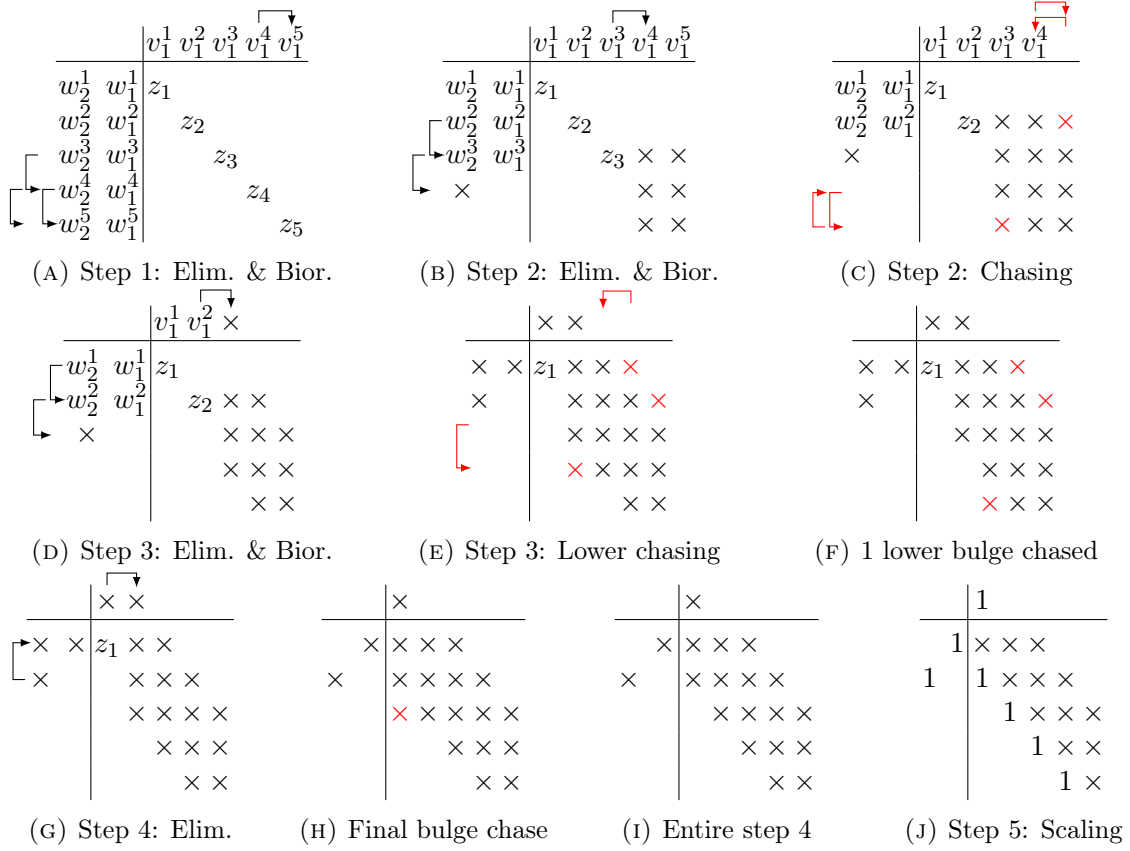

\textbf{Step 1.}
In the first step, the last element of $\bm{w}_2$ and $\bm{v}_1$, as well as the last two elements of $\bm{w}_1$ are eliminated.  More precisely, we compute eliminators $W^{-1}_{(2,5)}$, $W^{-1}_{(1,5)}$ and $W^{-1}_{(1,4)}$, where $W^{-1}_{(i,j)}$ is a lower eliminator that eliminates the $j\text{th}$ element of $\bm{w}_i$, i.e. the active part is a lower eliminator located in rows $j{-}1$ and $j$. Also, a matrix $V^{-1}_{(1,5)}$ is needed to eliminate the last element in $\bm{v}_1$. The total elimination matrices are $\widehat{W}_1^{-1} = W^{-1}_{(2,5)}W^{-1}_{(1,4)}W^{-1}_{(1,5)}$ and $\widehat{V}_1^{-1} = V^{-1}_{(1,5)}$. Applying these elimination matrices on the given vectors gives: 
\begin{align*}
	\widehat{W}_1^{-1} \begin{bmatrix} \bm{w}_1 & \bm{w}_2 \end{bmatrix}= \begin{bmatrix} \widehat{\bm{w}}_{1}&  \widehat{\bm{w}}_{2}\end{bmatrix}\qquad \text{and} \qquad \quad & \widehat{V}_1^{-1}\bm{v}_1 = \widehat{\bm{v}}_{1},
\end{align*}
where the last element of $\widehat{\bm{v}}_{1}$ and $ \widehat{\bm{w}}_{2}$ is zero, as well as the last two elements of $ \widehat{\bm{w}}_{1}$. The elimination is depicted in~\cref{fig:step0}, where the eliminators act on the vectors $\bm{w}_1$, $\bm{w}_2$, $\bm{v}_1$, and on the matrix $Z$. The eliminators on the left are applied to the left of $Z$, while those on the top are applied to the right. After the elimination step, we enforce biorthogonality as required by condition \ref{cond:2}. This is achieved by setting $W_1 = \widehat{W}_1 L_1^{\!-T}$ and $V_1 = \widehat{V}_1 U_1^{-1}$, where 
$\left({\widehat{W}_1}\right)^{\!T} \widehat{V}_1 = L_1 U_1$
is an LU decomposition.\footnote{We can not use a pivoted LU factorization. Hence, this a potential source for accuracy loss} With this choice, the biorthogonality is satisfied:
\begin{equation*}
    W_1^{\!T} V_1 = L_1^{-1} \widehat{W}_1^{T} \widehat{V}_1 U_1^{-1} = L_1^{-1} (L_1 U_1) U_1^{-1} = I.
\end{equation*}
Accordingly, the transformed $Z$ becomes
\begin{equation*}
    H^{(1)} = W_1^{\!T} Z V_1 = L_1^{-1} \widehat{W}_1^{T} Z \widehat{V}_1 U_1^{-1}.
\end{equation*}
It is essential that the biorthogonalization step preserves the zeros in $\bm{w}_1, \bm{w}_2, \bm{v}_1$ that were created in the previous steps. This is indeed the case as
\begin{align*}
	&W_1^{-1} \begin{bmatrix}  \bm{w}_1 & \bm{w}_2 \end{bmatrix} = L_1^{\!T} \left(\widehat{W}_1^{-1} \begin{bmatrix}  \bm{w}_1 & \bm{w}_2 \end{bmatrix}\right) =  L_1^{\!T} \begin{bmatrix} \widehat{\bm{w}}_{1}&  \widehat{\bm{w}}_{2}\end{bmatrix} = \begin{bmatrix}
		\bm{w}^{(1)}_{1} & \bm{w}^{(1)}_{2}
	\end{bmatrix},    \\
	&V_1^{-1} \bm{v}_1 = U_1 \left(\widehat{V}_1^{-1} \bm{v}_1\right) = U_1 \widehat{\bm{v}}_{\bm{1}} = \bm{v}^{(1)}_{1},
\end{align*} 
both retaining the created zeros since $L_1^{\!T}$ and $U_1$ are upper triangular.
As $Z_1$ has the desired bandwidth \cref{eq:Hessenberg-recurrence-matrix}, no additional elimination steps on $Z_1$, i.e., chasing, are necessary. The action of the operations on $\bm{w}_1$, $\bm{w}_2$, and $\bm{v}_1$, along with the related operation on $Z$, is illustrated in \cref{fig:step1}. Note that only the eliminators $\widehat{W}_1^{-1}$ and $\widehat{V}_1^{-1}$, indicated by arrow brackets, are shown. The matrices $L_1$ and $U_1$, used in the biorthogonalization, are omitted to avoid an overly complicated visual representation. Nevertheless, their action on $Z$ and the vectors $\bm{w}_1$, $\bm{w}_2$, and $\bm{v}_1$ is always illustrated.

\textbf{Step 2.}
The result after the elimination and biorthogonalization of step 2 is illustrated in \cref{fig:step2a}. The elimination differs slightly from step 1 as we only remove one element from $\bm{w}_{1}^{(1)}$ instead of two elements as in the previous step. More precisely, we apply  $\widehat{W}_2^{-1} = W^{-1}_{(2,4)}W^{-1}_{(1,3)}$ and $\widehat{V}_2^{-1} = V^{-1}_{(1,4)}$ such that 
\begin{equation*}
\widehat{W}_2^{-1} 
    \begin{bmatrix}
	\bm{w}_{1}^{(1)} & \bm{w}_{2}^{(1)}
\end{bmatrix} = \begin{bmatrix}
	\widehat{\bm{w}}_1^{(1)} & \widehat{\bm{w}}_2^{(1)}
\end{bmatrix}
\quad \text{and} \quad \widehat{V}_2^{-1}\bm{v}_{1}^{(1)} = \widehat{\bm{v}}_{1}^{(1)},
\end{equation*} where the last two elements of $\widehat{\bm{v}}_{1}^{(1)}$ and $\widehat{\bm{w}}_{2}^{(1)}$ are zero, as well as the last three elements of $\widehat{\bm{w}}_{1}^{(1)}$. Performing another LU decomposition gives the relations $ W_2 = \widehat{W}_2 L_2^{\!-T} $ and $ V_2 = \widehat{V}_2 U_2^{-1} $, ensuring that $ W_2^{\!T} V_2 = I $. Applying $W_2$ and $V_2$ to the vectors results in
\begin{align*}
	&W_2^{-1} \begin{bmatrix}  \bm{w}_1^{(1)} & \bm{w}_2^{(1)} \end{bmatrix} = L_2^{\!T} \left(\widehat{W}_2^{-1} \begin{bmatrix}  \bm{w}_1^{(1)} & \bm{w}_2^{(1)} \end{bmatrix}\right) =  L_2^{\!T} \begin{bmatrix} \widehat{\bm{w}}_1^{(1)} & \widehat{\bm{w}}_1^{(2)} \end{bmatrix} = \begin{bmatrix}
		\bm{w}^{(2)}_{1} & \bm{w}^{(2)}_{2}
	\end{bmatrix},    \\
	&V_2^{-1} \bm{v}_1^{(1)} = U_2 \left(\widehat{V}_2^{-1} \bm{v}_1^{(1)}\right) = U_2 \widehat{\bm{v}}_{1}^{(1)} = \bm{v}^{(2)}_{1},
\end{align*}
again retaining the desired zeros.
A second difference is that $\widetilde{H}^{(2)} = W_2^{T} H^{(1)} V_2$ has two unwanted elements which are indicated by bold crosses in \cref{fig:step2a}. We will refer to these unwanted elements as the \textit{bulges}. A chasing step has to be done to eliminate these bulges. The lower bulges, which are unwanted elements on the lower diagonal, are removed with an appropriate lower eliminator $W_{(2,1)}$ acting on the last two rows of $\widetilde{H}^{(2)}$. Since we must have a similarity transformation, we apply $W_{(2,1)}^{-1}$ on the right side of $\widetilde{H}^{(2)}$. A similar step with an upper eliminator $V_{(2,1)}$ acting on the last two columns of $\widetilde{H}^{(2)}$ is done to eliminate the element in the upper bulge. The transformation ${H}^{(2)} = V_{(2,1)}^{-1}W_{(2,1)}\widetilde{H}^{(2)}W_{(2,1)}^{\!-1} V_{(2,1)}$, shown in \cref{fig:step2b}, removes the unwanted elements in $\widetilde{H}^{(2)}$. Again, the eliminators in the chasing steps do not destroy the created zeros in $\bm{w}_{1}^{(2)}$ and $\bm{w}_{2}^{(2)}$. 

\textbf{Step 3.}
The elimination and biorthogonalization of step 3 are similar to the previous step, yielding $\widetilde{H}^{(3)}= W_3^TH^{(2)}V_3$. In \cref{fig:step3a}, one can see that the bulge now consists of 2 upper elements and 1 lower element. The chasing procedure now becomes a bit more technical. The chasing step is illustrated for the lower bulges, as the process for the upper bulges is analogous. As depicted with the red arrow, one can eliminate the lower bulges by a lower eliminator $W_{3,1}$ acting on the left of $\widetilde{H}^{(4)}$. To enforce biorthogonality, we apply the inverse on the right, resulting in $W_{(3,1)}\widetilde{H}^{(4)}W_{(3,1)}^{-1}$ which is shown in \cref{fig:step3b}. Now, the biorthogonality constraint introduced a new lower bulge, positioned in the last row. Because of this bulge movement, this step is commonly referred to as \emph{bulge chasing} (see, e.g., \cite[Chapter~7]{GolVanL13}). The lower-positioned bulge can be removed by a transformation similar to the chasing part of step 2. \Cref{fig:step4a} shows the situation after the entire chasing procedure. In general, in the chasing step of step $i$ where $1<i<N-1$, we perform $i{-}1$ eliminations in both the lower and upper part of the matrix.

\textbf{Step 4.}
The penultimate step is special, as we use upper eliminators to create the necessary zeros, instead of relying on lower eliminators. More specifically, take $W^{-1}_4$ as an upper eliminator such that $W^{-1}_4\bm{w}^{(3)}_{2} = c \,\bm{e}_2$, where $c$ is a nonzero constant; i.e., the first element of $\bm{w}^{(3)}_{2}$ is eliminated. We also need to eliminate the second element in $v_1^{(3)}$ in this step; see \cref{fig:step4a}. It turns out that $V_4^{-1} = W^{T}_4$ can be used for this, which at the same time ensures the biorthogonality condition for this step. To see why this choice is justified, we note that during all steps, the orthogonality of $\v_1^{(k)}$ and $\w_2^{(k)}$ is preserved. In particular, 
\begin{equation*}
    0 = \left({\bm{w}^{(4)}_{2}}\right)^T {\bm{v}^{(4)}_{1}} = \left(W_4^{-1}{\bm{w}^{(3)}_{2}}\right)^T \left(W_4^T{\bm{v}^{(3)}_{1}}\right) = c \bm{e}_2^T \left(W_4^T{\bm{v}^{(3)}_{1}}\right).
\end{equation*}
Thus, the second element of $W_4^T{\bm{v}^{(3)}_{1}}$ is eliminated as $c$ is not zero. From \cref{fig:step4b}, it can be seen that $\widetilde{H}^{(4)}= W_4^TH^{(3)}V_4$ still has one lower bulge element. A final lower bulge-chasing step is required to obtain the configuration shown in \cref{fig:step4c}.

\textbf{Step 5.}
The last part of the algorithm applies appropriate scaling matrices. We apply matrices $S_1,S_2,S_3$ such that $S_1 \bm{w}^{(4)}_1 = \bm{e}_1$, $S_2 \bm{w}^{(4)}_2 = \bm{e}_2$ and $S_3 \bm{v}^{(4)}_1 = \bm{e}_1$; see \cref{fig:step5}. Recall from the initial conditions that we have $\bm{w}_2^T Z \bm{v}_1 = 1$. At this stage of the algorithm, we satisfy condition \ref{cond:1}, and thus 
\begin{equation*}
    \bm{e}^T_2H^{(4)}\bm{e}_1 = \bm{e}^T_2W_{5}^TZV_{5}\bm{e}_1 = \bm{w}_2^T Z \bm{v}_1 = 1. 
\end{equation*}
The second scaling step applies \rev{the diagonal} similarity transformation \rev{described in \cref{sec:Krylov}} on $H^{(5)}$ such that all elements on the subdiagonal are one. To this end, we compute $H^{(5)}=D_sH^{(4)}D_s^{-1}$ with ($H_{ij}^{(k)}$ denotes the element $(i,j)$ of $H^{(k)}$):
\begin{equation}
\label{eq:diagonal-scaling-small}
    D_s = \diag \left(1,1,\frac{1}{H^{(4)}_{3,2}}, \frac{1}{H^{(4)}_{3,2}H^{(4)}_{4,3}}, \frac{1}{H^{(4)}_{3,2}H^{(4)}_{4,3}H^{(4)}_{5,4}}\right).
\end{equation}
\begin{remark}
    We \rev{choose} to not apply this scaling to the basis vectors $W$ and $V$, as it is an ill-conditioned transformation due to the typically large magnitude variation of the entries in $D_s$\rev{; see \cref{sec:Krylov} for more details.}
\end{remark}

\section{Numerical experiments}\label{sec:numericalexperiments}
In this section, we compare our numerical algorithms with the explicit formulae given by Arvesú, Coussement, and Van Assche~\cite{ArvCousVanAssche2003} for the recurrence matrix associated with two examples of discrete multiple orthogonal polynomials, namely the Kravchuk and Hahn MOPs. In \cref{sec:conditioning/stability},
we illustrate --- using the multi-precision computing toolbox Advanpix in MATLAB --- that the IEP linked to these Kravchuk and Hahn polynomials is ill-conditioned. Furthermore, we also investigate the stability of our algorithms. Finally, we consider IEPs which are better conditioned than 
the ones related to
Kravchuk and Hahn MOPs to illustrate that our algorithms perform well. All MATLAB codes are publicly \rev{available}.\footnote{\url{https://gitlab.kuleuven.be/numa/public/IEP-MOP}}

\subsection{Kravchuk and Hahn Multiple Orthogonal Polynomials}\label{subsec:HahnandKravchuk}
Arvesú et al.~\cite{ArvCousVanAssche2003} gave explicit expressions for certain discrete multiple orthogonal polynomials and their recurrence coefficients. These examples include Charlier, Meixner, Kravchuk and Hahn MOPs. We tested Hahn and Kravchuk MOPs, as our setting assumes finitely many nodes in the discrete measures.  The measures corresponding to these MOPs are examples of AT-systems for the case $r=2$ (cf.~\cite{ArvCousVanAssche2003}).
	\subsubsection*{Example 1: Multiple Kravchuk polynomials}
		Kravchuk polynomials arise when the discrete measures are binomial distributions on the integers. More concretely,
		\begin{equation*}
			\mu_j = \sum_{i=0}^{N-1} p_j^i(1-p_j)^{N-i}\delta_{i}, \quad j\in\{1,2\},
		\end{equation*}
	where $0<p_j<1$ are all distinct. We use $p_1 = 0.4$ and $p_2=0.5$ in all our experiments. Explicit formulae for the type II multiple Kravchuk polynomials $K_{n1,n2}^{p1,p2,N}(x)$ can be derived. This is done by defining a suitable \textit{raising operator} and corresponding \textit{Rodrigues formula}~\cite[Sec.~4.4]{ArvCousVanAssche2003}. For the multi-index $(n_{1},n_{2})$, one obtains
	\begin{multline*}
		K_{n1,n2}^{p1,p2,N}(x)= \\
         p_1^{n_1} p_2^{n_2} (-N-1)_{n_1+n_2} \sum_{j=0}^{n_1+n_2} \sum_{k=0}^{j} 
		\frac{(-n_1)_k}{k!} 
		\left( \frac{1}{p_1} \right)^k 
		\frac{(-n_2)_{j-k}}{(j-k)!} 
		\left( \frac{1}{p_2} \right)^{j-k} 
		\frac{(-x)_j}{(-N-1)_j},
	\end{multline*}
    where $ (c)_j $ denotes the Pochhammer function of $c$ , given by $ (c)_j = \prod_{i=0}^{j-1} (c + i) $ for $ j > 0 $, and $ (c)_0 = 1 $. 
	\subsubsection*{Example 2: Multiple Hahn polynomials}
	We obtain Hahn polynomials when we take hypergeometric distributions on the integers as measures. We get
	\begin{equation*}
		\mu_j = \sum_{i=0}^{N-1} \frac{(\beta_j+1)_i}{i!}\frac{(\gamma+1)_{N-i-1}}{(N-1-i)!} \delta_i, \quad \beta_j >-1,\quad \gamma > -1,\quad  j\in\{1,2\},
	\end{equation*}
	with all $\beta_j$ distinct. We use $\beta_1=1$, $\beta_2=1.5$ and $\gamma=1$ in all our experiments. Also, for this AT-system, one can find explicit formulae~\cite{ArvCousVanAssche2003}. 
    
    We compare our algorithms by computing the forward relative error 
\begin{equation}\label{eq:relerror}
    e_N = \frac{\left \lVert H_N - \widehat{H}_N \right \rVert_2}{\left \lVert H_N \right \rVert_2},
\end{equation}
where $\widehat{H}_N$ is the numerically obtained solution and $H_N$ is the exact solution~\cite{ArvCousVanAssche2003}. In \cref{fig:experiment1}, we show the forward relative error $e_N$ and the loss of biorthogonality $\left \lVert I_N - W^TV \right \rVert_2$, plotted over dimensions $N = 5,\ldots,30$ for the Kravchuk and Hahn MOPs. From this, we can draw a couple of observations. The forward relative error $e_N$ has similar growth for \texttt{IEP\_CORE} and \texttt{IEP\_KRYLREORTH(full)}, while the error is larger for \texttt{IEP\_KRYL}. As explained in previous sections, the bases $W$ and $V$ for \texttt{IEP\_KRYL} correspond to monic type II MOPs which are ill-conditioned. In contrast, we do not scale $W$ and $V$ for \texttt{IEP\_CORE} and \texttt{IEP\_KRYLREORTH(full)} in order to have a better conditioned basis, and thus they are diagonal scalings of the basis vectors connected to monic type II MOPs. The loss of biorthogonality for \texttt{IEP\_CORE} and \texttt{IEP\_KRYLREORTH(full)} remains small with increasing dimension, while it increases exponentially for the ill-conditioned basis vectors in \texttt{IEP\_KRYL}. In \cref{fig:experiment6}, the \rev{relative entry-wise} distance between $H_{20}$ and the numerical solution $\widehat{H}_{20}$ for the different diagonals, computed with \texttt{IEP\_KRYLREORTH(full)}, is visualized. Although the initial columns of ${H}_N$ are approximated with high precision, the error increases exponentially as the column index grows. The problems appear to be highly ill-conditioned as will be analyzed in more detail in the next section.
\begin{figure}
    \centering
    \tikzsetnextfilename{experiment1}
\pgfplotsset{height=0.48\linewidth,width=0.49\linewidth}
\pgfplotsset{major grid style={dotted,gray}}

\noindent%
\begin{tikzpicture}[scale=1,baseline]%
\begin{groupplot}[
    group style={
        group size=2 by 1,
        horizontal sep=2cm,
    },
]

\nextgroupplot[label style={font=\scriptsize}, tick label style={font=\tiny}, legend columns = 1,mark options={scale=0.6}, legend style={at = {(0.01,0.99)}, anchor = north west,row sep=-0.2pt,nodes={scale=0.485, transform shape}},legend cell align={left}, grid=major, xlabel={$N$}, xmin = 5, xmax = 30, ymin=1e-17,ymax=10e15, ymode=log, xtick distance=5, ytick distance = 10000, minor ytick={1e-16,1e-15,1e-14,1e-12,1e-11,1e-10,1e-8,1e-7,1e-6,1e-4,1e-3,1e-2,1e0,1e1,1e2,1e4,1e5,1e6,1e8,1e9,1e10,1e12,1e13,1e14}]

\addplot[red, mark=*, solid] table[x={N}, y={err_H_core}] {Experiments/experiment1/data_experiment1_kravchuk.dat};
\addlegendentry{\small $e_N$ \texttt{IEP\_CORE}}
\addplot[blue, mark=diamond*, solid] table[x={N}, y={err_H_krylov}] {Experiments/experiment1/data_experiment1_kravchuk.dat};
\addlegendentry{\small $e_N$ \texttt{IEP\_KRYL}}
\addplot[green!50!black, mark=triangle*, solid] table[x={N}, y={err_H_krylov_reorth}] {Experiments/experiment1/data_experiment1_kravchuk.dat};
\addlegendentry{\small $e_N$ \texttt{IEP\_KRYLREORTH(full)}}
\addplot[red, dashed] table[x={N}, y={err_bior_core}] {Experiments/experiment1/data_experiment1_kravchuk.dat};
\addlegendentry{\small $\|W^T V - I_n\|_2$ \texttt{IEP\_CORE}}
\addplot[blue, dashed] table[x={N}, y={err_bior_krylov}] {Experiments/experiment1/data_experiment1_kravchuk.dat};
\addlegendentry{\small $\|W^T V - I_n\|_2$ \texttt{IEP\_KRYL}}
\addplot[green!70!black, dashed] table[x={N}, y={err_bior_krylov_reorth}] {Experiments/experiment1/data_experiment1_kravchuk.dat};
\addlegendentry{\small $\|W^T V - I_n\|_2$ \texttt{IEP\_KRYLREORTH(full)}}
\addplot[black, solid] table[x={N}, y={error_cond}] {Experiments/experiment1/data_experiment1_kravchuk.dat};
\addlegendentry{\small Conditioning error}

\nextgroupplot[label style={font=\scriptsize}, tick label style={font=\tiny}, mark options={scale=0.6},grid=major, xlabel={$N$}, xmin = 5, xmax = 30, ymin=10e-18,ymax=10e15, ymode=log, xtick distance=5, ytick distance = 10000,minor ytick={1e-16,1e-15,1e-14,1e-12,1e-11,1e-10,1e-8,1e-7,1e-6,1e-4,1e-3,1e-2,1e0,1e1,1e2,1e4,1e5,1e6,1e8,1e9,1e10,1e12,1e13,1e14},]

\addplot[red, mark=*, solid] table[x={N}, y={err_H_core}] {Experiments/experiment1/data_experiment1_hahn.dat};
\addplot[blue, mark=diamond*, solid] table[x={N}, y={err_H_krylov}] {Experiments/experiment1/data_experiment1_hahn.dat};
\addplot[green!50!black, mark=triangle*, solid] table[x={N}, y={err_H_krylov_reorth}] {Experiments/experiment1/data_experiment1_hahn.dat};
\addplot[red, dashed] table[x={N}, y={err_bior_core}] {Experiments/experiment1/data_experiment1_hahn.dat};
\addplot[blue, dashed] table[x={N}, y={err_bior_krylov}] {Experiments/experiment1/data_experiment1_hahn.dat};
\addplot[green!70!black, dashed] table[x={N}, y={err_bior_krylov_reorth}] {Experiments/experiment1/data_experiment1_hahn.dat};
\addplot[black, solid] table[x={N}, y={error_cond}] {Experiments/experiment1/data_experiment1_hahn.dat};

\end{groupplot}
\end{tikzpicture}
    \caption{Accuracy and biorthogonalization loss of the various algorithms for the Kravchuk (left) and Hahn (right) MOPs.}
    \label{fig:experiment1}
\end{figure}
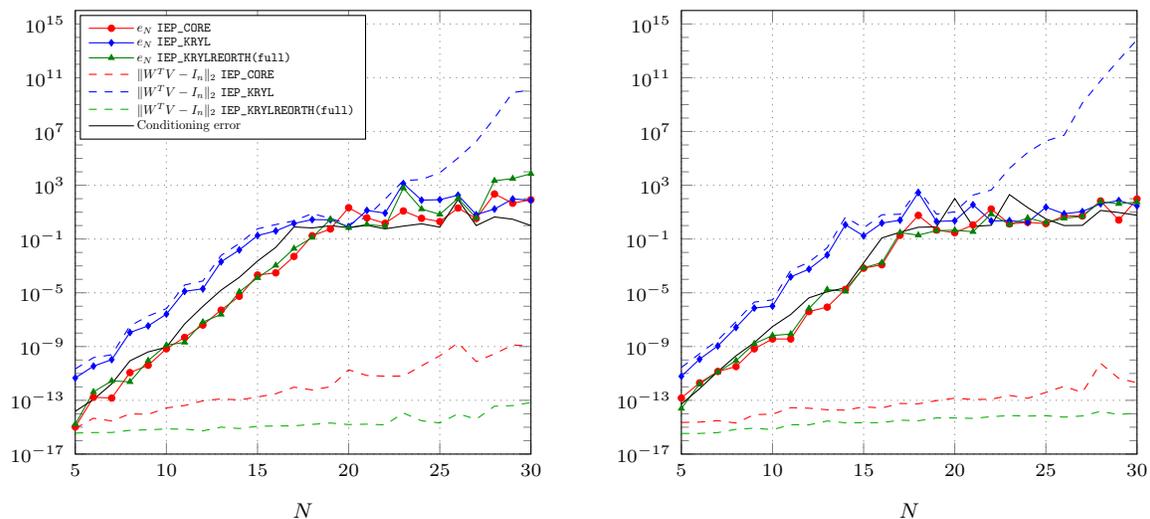
\begin{figure}
    \centering
    \tikzsetnextfilename{experiment6}
\pgfplotsset{height=0.47\linewidth,width=0.47\linewidth}
\pgfplotsset{major grid style={dotted,gray}}

\noindent%
\begin{tikzpicture}[scale=1,baseline]%
\begin{groupplot}[
    group style={
        group size=2 by 1,
        horizontal sep=2cm,
    },
]

\nextgroupplot[label style={font=\tiny}, tick label style={font=\tiny},mark options={scale=0.6}, legend columns = 1, legend style={at = {(0.01,0.99)}, anchor = north west,row sep=-0.2pt,nodes={scale=0.55, transform shape}},legend cell align={left}, grid=major, xlabel={Column index $j$}, xmin = 1, xmax = 20, ymin=1e-17,ymax=1e3, ymode=log, xtick distance=5, ytick distance = 10000, minor ytick={1e-16, 1e-15, 1e-14, 1e-12, 1e-11, 1e-10, 1e-8, 1e-7, 1e-6, 1e-4, 1e-3, 1e-2, 1e0, 1e1, 1e2, 1e4, 1e5, 1e6, 1e8, 1e9, 1e10, 1e12, 1e13, 1e14, 1e16, 1e17, 1e18, 1e20, 1e21, 1e22, 1e24, 1e25, 1e26,1e28, 1e29, 1e30,1e32, 1e33, 1e34}]

\addplot[blue, mark=*, solid] table[x={N}, y={err_diag_kravchuk}] {Experiments/experiment6/data_experiment6_k.dat};
\addlegendentry{\small Main diagonal}
\addplot[red, mark=diamond*, solid] table[x={N}, y={err_diag1_kravchuk}] {Experiments/experiment6/data_experiment6_k1.dat};
\addlegendentry{\small 1\textsuperscript{st} superdiagonal}
\addplot[green!50!black, mark=triangle*, solid] table[x={N}, y={err_diag2_kravchuk}]{Experiments/experiment6/data_experiment6_k2.dat};
\addlegendentry{\small 2\textsuperscript{nd} superdiagonal}

\nextgroupplot[label style={font=\tiny}, tick label style={font=\tiny}, mark options={scale=0.6},legend columns = 1,legend style={at = {(0.01,0.99)}, anchor = north west,row sep=-0.2pt,nodes={scale=0.55, transform shape}},legend cell align={left},grid=major, xlabel={Column index $j$}, xmin = 1, xmax = 20, ymin=1e-17,ymax=1e3, ymode=log, xtick distance=5, ytick distance = 10000,minor ytick={1e-16, 1e-15, 1e-14, 1e-12, 1e-11, 1e-10, 1e-8, 1e-7, 1e-6, 1e-4, 1e-3, 1e-2, 1e0, 1e1, 1e2, 1e4, 1e5, 1e6, 1e8, 1e9, 1e10, 1e12, 1e13, 1e14, 1e16, 1e17, 1e18, 1e20, 1e21, 1e22, 1e24, 1e25, 1e26,1e28, 1e29, 1e30,1e32, 1e33, 1e34}]

\addplot[blue, mark=*, solid] table[x={N}, y={err_diag_hahn}] {Experiments/experiment6/data_experiment6_h.dat};
\addplot[red, mark=diamond*, solid] table[x={N}, y={err_diag1_hahn}] {Experiments/experiment6/data_experiment6_h1.dat};
\addplot[green!50!black, mark=triangle*, solid] table[x={N}, y={err_diag2_hahn}]{Experiments/experiment6/data_experiment6_h2.dat};

\end{groupplot}
\end{tikzpicture}
    \caption{\rev{Relative entry-wise accuracy $\left\lvert\frac{[H_{20}- \widehat{H}_{20}]_{i,j}}{[H_{20}]_{i,j}}\right\rvert$} of \texttt{IEP\_KRYLREORTH(full)} for Kravchuk MOPs (left) and Hahn MOPs (right).}
    \label{fig:experiment6}
\end{figure}

\subsection{Conditioning and stability}\label{sec:conditioning/stability}
The errors in \cref{fig:experiment1} show that for values of $N$ between 15 and 20, all decimal digits of precision for the entries in $\widehat{H}_N$ are lost. This loss may be attributed either to instability of the methods or to ill-conditioning of the underlying problem.
To investigate this, we examine the conditioning of the problem where the input consists of the nodes $z_i$ and the weights $\alpha_{j,i}$ of the two discrete measures, and the output is the banded Hessenberg matrix $H_N$. Given this input, we compute the matrix $H_N$ in quadruple precision using the Advanpix Multiprecision Computing Toolbox.\footnote{\url{https://www.advanpix.com/}} 
\rev{Here and in the remainder of the text, we approximate exact
arithmetic with computations in quadruple precision.}
Next, we perturb the input data with relative errors on the order of \rev{machine epsilon for double precision} and compute the solution of the problem $\widetilde{H}_N$ in quadruple precision\rev{, and treat it as the exact solution}. The relative difference between $ H_N $ and $ \widetilde{H}_N $, shown as the black line in \cref{fig:experiment1}, represents the conditioning error. The forward relative errors $e_N$ are of the same order of magnitude as the conditioning error, indicating that all algorithms are weakly stable for this problem.
 
\begin{figure}
    \centering
    \tikzsetnextfilename{experiment2}
\pgfplotsset{height=0.47\linewidth,width=0.48\linewidth}
\pgfplotsset{major grid style={dotted,gray}}

\noindent%
\begin{tikzpicture}[scale=1,baseline]%
\begin{groupplot}[
    group style={
        group size=2 by 1,
        horizontal sep=2cm,
    },
]

\nextgroupplot[label style={font=\tiny}, tick label style={font=\tiny}, legend columns = 1,mark options={scale=0.6}, legend style={at = {(0.4,0.83)}, anchor = north west,row sep=-0.2pt,nodes={scale=0.485, transform shape}},legend cell align={left}, grid=major, xlabel={$N$}, ylabel={$\|\bm z - \widehat{\bm z}\|_2 / \|\bm z\|_2$}, xmin = 5, xmax = 30, ymin=1e-17,ymax=1e0, ymode=log, xtick distance=5, ytick distance = 10000, minor ytick = {1e-16, 1e-15, 1e-14, 1e-12, 1e-11, 1e-10, 1e-8, 1e-7, 1e-6, 1e-4, 1e-3, 1e-2, 1e0, 1e1, 1e2, 1e4, 1e5, 1e6, 1e8, 1e9, 1e10, 1e12, 1e13, 1e14, 1e16, 1e17, 1e18, 1e20, 1e21, 1e22, 1e24, 1e25, 1e26,1e28, 1e29, 1e30,1e32, 1e33, 1e34}]

\addplot[red, mark=*, solid] table[x={N}, y={err_nodes_core}] {Experiments/experiment2/data_experiment2.dat};
\addlegendentry{\small \texttt{IEP\_CORE}}
\addplot[blue, mark=diamond*, solid] table[x={N}, y={err_nodes_kryl}] {Experiments/experiment2/data_experiment2.dat};
\addlegendentry{\small \texttt{IEP\_KRYL}}
\addplot[green!50!black, mark=triangle*, solid] table[x={N}, y={err_nodes_kryl_reorth}] {Experiments/experiment2/data_experiment2.dat};
\addlegendentry{\small \texttt{IEP\_KRYLREORTH(full)}}

\nextgroupplot[label style={font=\tiny}, tick label style={font=\tiny}, legend columns = 1, mark options={scale=0.6}, legend style={at = {(0.37,0.83)}, anchor = north west,row sep=-0.2pt,nodes={scale=0.485, transform shape}},legend cell align={left},grid=major, xlabel={$N$},ylabel={$\|\bm{\alpha}_j - \widehat{\bm{\alpha}}_j\|_2 / \|\boldsymbol{\alpha}_j\|_2$}, xmin = 5, xmax = 30, ymin=10e-18,ymax=1e0, ymode=log, xtick distance=5, ytick distance = 10000,minor ytick={1e-16, 1e-15, 1e-14, 1e-12, 1e-11, 1e-10, 1e-8, 1e-7, 1e-6, 1e-4, 1e-3, 1e-2, 1e0, 1e1, 1e2, 1e4, 1e5, 1e6, 1e8, 1e9, 1e10, 1e12, 1e13, 1e14, 1e16, 1e17, 1e18, 1e20, 1e21, 1e22, 1e24, 1e25, 1e26,1e28, 1e29, 1e30,1e32, 1e33, 1e34}]

\addplot[red, mark=*, solid] table[x={N}, y={err_w1_core}] {Experiments/experiment2/data_experiment2.dat};
\addlegendentry{\small  $j=1$ \texttt{IEP\_CORE}}
\addplot[blue, mark=diamond*, solid] table[x={N}, y={err_w1_kryl}] {Experiments/experiment2/data_experiment2.dat};
\addlegendentry{\small $j=1$ \texttt{IEP\_KRYL}}
\addplot[green!50!black, mark=triangle*, solid] table[x={N}, y={err_w1_kryl_reorth}] {Experiments/experiment2/data_experiment2.dat};
\addlegendentry{\small $j=1$ \texttt{IEP\_KRYLREORTH(full)}}

\addplot[red, mark=*, dashed] table[x={N}, y={err_w2_core}] {Experiments/experiment2/data_experiment2.dat};
\addlegendentry{\small $j=2$ \texttt{IEP\_CORE}}
\addplot[blue, mark=diamond*, dashed] table[x={N}, y={err_w2_kryl}] {Experiments/experiment2/data_experiment2.dat};
\addlegendentry{\small $j=2$ \texttt{IEP\_KRYL}}
\addplot[green!50!black, mark=triangle*, dashed] table[x={N}, y={err_w2_kryl_reorth}] {Experiments/experiment2/data_experiment2.dat};
\addlegendentry{\small $j=2$ \texttt{IEP\_KRYLREORTH(full)}}

\end{groupplot}
\end{tikzpicture}
    \caption{Comparison of backward error on nodes (left) and weights (right) for multiple Hahn polynomials. Note that the dashed lines in the right figure almost coincide with the solid lines.}
    \label{fig:experiment2}
\end{figure}
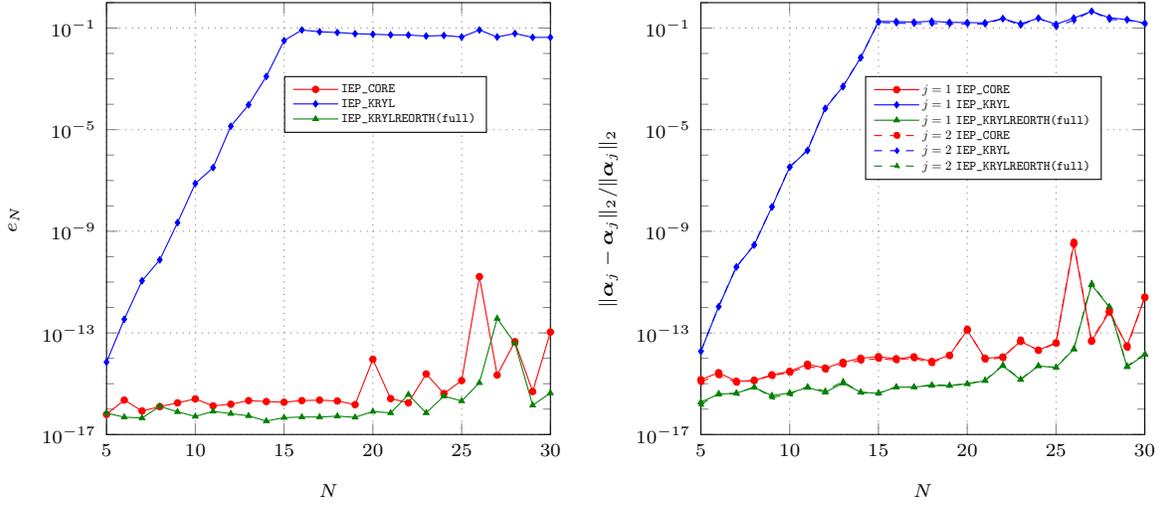

To assess the backward stability we check if the computed solution $\widehat{H}_N$
corresponds to the exact solution of slightly perturbed input data 
$\hat{z}_i$ and $\hat{\alpha}_{j,i}$ where the
relative perturbation is of the size of the machine precision. \rev{More specifically, we compute the eigenvalues and eigenvectors of $\widehat{H}_N$ in quadruple precision, obtaining the nodes $\widehat{\bm z} = [\hat{z}_1, \dots, \hat{z}_N]$ as the eigenvalues. The weight vectors $\widehat{\bm \alpha}_1$ and $\widehat{\bm \alpha}_2$ are then computed using \cref{theorem:VanAsscheThm2}. 
In \cref{fig:experiment2}, we show the relative backward errors
\begin{equation}
\label{eq:backward-errors}
\frac{\|\bm z - \widehat{\bm z}\|_2}{\|\bm z\|_2},
\quad
\frac{\|\bm \alpha_j - \widehat{\bm \alpha}_j\|_2}{\|\bm \alpha_j\|_2},
\end{equation}
where $\bm z$ and $\bm \alpha_j$ contain the original nodes and weights. 
 We use the input data from} the multiple Hahn
 polynomials with the same parameters as before. The results presented in the figure indicate that \texttt{IEP\_CORE} and \texttt{IEP\_KRYLREORTH(full)} exhibit backward stability for most values of $N<25$.
For $N \geq 25$, backward instability is observed for both methods. In contrast, \texttt{IEP\_KRYL} is not backward stable for any of the considered values of $N$.

\subsection{Other IEPs}
In this section, we solve other IEPs which are better conditioned than the Kravchuk and Hahn examples. 
\rev{We solve the IEP associated with monic type II MOPs, i.e., with ones on the subdiagonal of the recurrence matrix, which solution is unique as stated in \cref{prop:uniqueness-IEP}.} 
Although we don't have the explicit formulae for the recurrence coefficients of these examples, we approximate the exact solution using quadruple precision in Advanpix. This approximated exact recurrence matrix is denoted by $H_N$ and the numerical approximation using double precision with $\widehat{H}_N$. Similar as in \cref{subsec:HahnandKravchuk}, we compare the algorithms by using the forward relative error $e_N$ \cref{eq:relerror} and loss of biorthogonality. Note that the examples below are not guaranteed to be AT-systems. Nevertheless, we numerically verified the rank of the moment matrix (see equation \cref{eq:momentmatrix}) to ensure a normal index, which guarantees the existence of a unique set of MOPs.

\begin{itemize}
    \item \textbf{Equidistant nodes in $[-1,1]$ and random weights $\boldsymbol{\alpha}_j \sim \mathcal{U}(1, 2)^N$:}
    
$\mathcal{U}(a,b)$ denotes the uniform distribution on the interval $(a,b)$. \rev{For each value of $N$, we take the previously generated weights $\bm\alpha_j\in \mathcal{U}(1,2)^{N-1}$ and append one additional weight sampled from $\mathcal{U}(1,2)$. Since the weights are random, all reported results are averaged over 20 independent runs. } The results in \cref{fig:experiment3} indicate that the reorthogonalized version \texttt{IEP\_KRYLREORTH(full)} outperforms \rev{\texttt{IEP\_KRYL}}. The error associated with \texttt{IEP\_CORE} appears to be highly sensitive to the randomness in the weights. 
\rev{We show the loss of orthogonality before and after the diagonal scaling which is described in equation \cref{eq:diagonal-scaling-ones}. Although the recurrence matrix is well approximated after the scaling in \texttt{IEP\_KRYLREORTH(full)}, the corresponding scaled bases show a significant loss of orthogonality compared to the unscaled bases, with the same growth rate as for \texttt{IEP\_KRYL}.}

\begin{figure}
    \centering
    \tikzsetnextfilename{experiment3}
\pgfplotsset{height=0.47\linewidth,width=0.48\linewidth}
\pgfplotsset{major grid style={dotted,gray}}

\noindent%
\begin{tikzpicture}[scale=1,baseline]%
\begin{groupplot}[
    group style={
        group size=2 by 1,
        horizontal sep=2cm,
    },
]

\nextgroupplot[label style={font=\tiny}, tick label style={font=\tiny},mark options={scale=0.6}, legend columns = 1, legend style={at = {(0.01,0.99)}, anchor = north west,row sep=-0.2pt,nodes={scale=0.55, transform shape}},legend cell align={left}, grid=major, xlabel={$N$},ylabel={$e_N$}, xmin = 5, xmax = 50, ymin=1e-17,ymax=1e7, ymode=log, xtick distance=10, ytick distance = 10000, minor ytick={1e-16, 1e-15, 1e-14, 1e-12, 1e-11, 1e-10, 1e-8, 1e-7, 1e-6, 1e-4, 1e-3, 1e-2, 1e0, 1e1, 1e2, 1e4, 1e5, 1e6, 1e8, 1e9, 1e10, 1e12, 1e13, 1e14, 1e16, 1e17, 1e18, 1e20, 1e21, 1e22, 1e24, 1e25, 1e26,1e28, 1e29, 1e30,1e32, 1e33, 1e34}]

\addplot[red, mark=*, solid] table[x={N}, y={err_core}] {Experiments/experiment3/errors_final.dat};
\addlegendentry{\small \texttt{IEP\_CORE}}
\addplot[blue, mark=diamond*, solid] table[x={N}, y={err_kryl}] {Experiments/experiment3/errors_final.dat};
\addlegendentry{\small \texttt{IEP\_KRYL}}
\addplot[green!50!black, mark=triangle*, solid] table[x={N}, y={err_kryl_reorth}] {Experiments/experiment3/errors_final.dat};
\addlegendentry{\small  \texttt{IEP\_KRYLREORTH(full)}}

\nextgroupplot[label style={font=\tiny}, tick label style={font=\tiny}, mark options={scale=0.6},legend columns = 1,legend style={at = {(0.01,0.99)}, anchor = north west,row sep=-0.2pt,nodes={scale=0.55, transform shape}},legend cell align={left},grid=major, xlabel={$N$},ylabel={$\|W^TV-I_N\|_2$}, xmin = 5, xmax = 50, ymin=1e-17,ymax=1e20, ymode=log, xtick distance=10, ytick distance = 10000,minor ytick={1e-16, 1e-15, 1e-14, 1e-12, 1e-11, 1e-10, 1e-8, 1e-7, 1e-6, 1e-4, 1e-3, 1e-2, 1e0, 1e1, 1e2, 1e4, 1e5, 1e6, 1e8, 1e9, 1e10, 1e12, 1e13, 1e14, 1e16, 1e17, 1e18, 1e20, 1e21, 1e22, 1e24, 1e25, 1e26,1e28, 1e29, 1e30,1e32, 1e33, 1e34}]

\addplot[red, mark=*, solid] table[x={N}, y={err_bior_core}] {Experiments/experiment3/errors_final.dat};
\addlegendentry{\small  \texttt{IEP\_CORE}}
\addplot[blue, mark=diamond*, solid] table[x={N}, y={err_bior_kryl}] {Experiments/experiment3/errors_final.dat};
\addlegendentry{\small \texttt{IEP\_KRYL}}
\addplot[green!50!black, mark=triangle*, solid] table[x={N}, y={err_bior_kryl_reorth}] {Experiments/experiment3/errors_final.dat};
\addlegendentry{\small  \texttt{IEP\_KRYLREORTH(full)}}
\addplot[red, mark=diamond*, dashed] table[x={N}, y={err_bior_core_scaled}] {Experiments/experiment3/errors_final.dat};
\addlegendentry{\small \texttt{IEP\_CORE} \textbf{scaled}}
\addplot[green!50!black, mark=triangle*, dashed] table[x={N}, y={err_bior_kryl_reorth_scaled}] {Experiments/experiment3/errors_final.dat};
\addlegendentry{\small  \texttt{IEP\_KRYLREORTH(full) \textbf{scaled}}}

\end{groupplot}
\end{tikzpicture}
    \caption{Error comparison (left) and loss of biorthogonality (right) for equidistant nodes in $[-1,1]$ and random weights, averaged over 20 runs.}
    \label{fig:experiment3}
\end{figure}
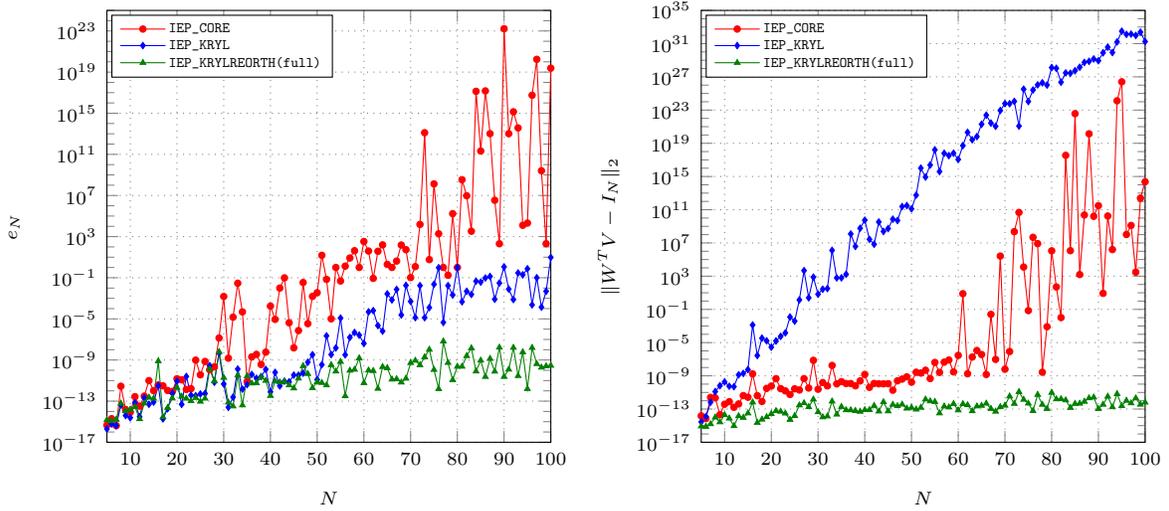

\item \textbf{Chebyshev nodes in $[-1,1]$ and random weights $\boldsymbol{\alpha}_j \sim \mathcal{U}(1, 2)^N$:}

Due to the good approximation properties of Chebyshev nodes, one would expect that this gives a better-conditioned inverse eigenvalue problem. This explains why in \cref{fig:experiment4}, the errors $e_N$ of \texttt{IEP\_KRYL} and \texttt{IEP\_KRYLREORTH(full)} increase less rapidly compared to the errors for the equidistant nodes. Both Krylov methods still clearly outperform \texttt{IEP\_CORE}.

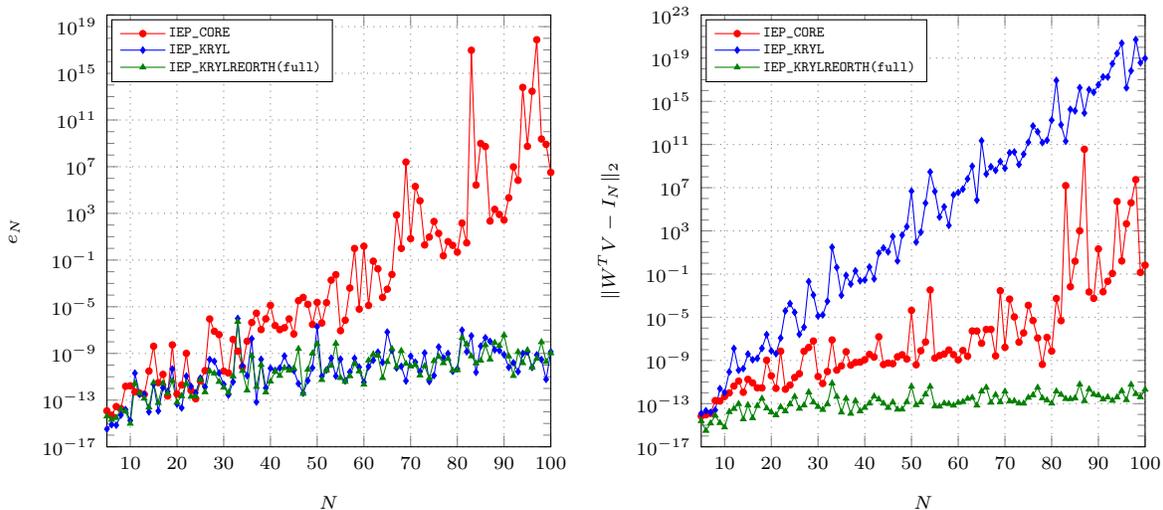
\begin{figure}
    \centering
    \tikzsetnextfilename{experiment4}
\pgfplotsset{height=0.47\linewidth,width=0.48\linewidth}
\pgfplotsset{major grid style={dotted,gray}}

\noindent%
\begin{tikzpicture}[scale=1,baseline]%
\begin{groupplot}[
    group style={
        group size=2 by 1,
        horizontal sep=2cm,
    },
]

\nextgroupplot[label style={font=\tiny}, tick label style={font=\tiny},mark options={scale=0.6}, legend columns = 1, legend style={at = {(0.01,0.99)}, anchor = north west,row sep=-0.2pt,nodes={scale=0.55, transform shape}},legend cell align={left}, grid=major, xlabel={$N$},ylabel={$e_N$}, xmin = 5, xmax = 50, ymin=1e-17,ymax=1e2, ymode=log, xtick distance=10, ytick distance = 10000, minor ytick={1e-16, 1e-15, 1e-14, 1e-12, 1e-11, 1e-10, 1e-8, 1e-7, 1e-6, 1e-4, 1e-3, 1e-2, 1e0, 1e1, 1e2, 1e4, 1e5, 1e6, 1e8, 1e9, 1e10, 1e12, 1e13, 1e14, 1e16, 1e17, 1e18, 1e20, 1e21, 1e22, 1e24, 1e25, 1e26,1e28, 1e29, 1e30,1e32, 1e33, 1e34}]

\addplot[red, mark=*, solid] table[x={N}, y={err_core}] {Experiments/experiment4/errors_final.dat};
\addlegendentry{\small \texttt{IEP\_CORE}}
\addplot[blue, mark=diamond*, solid] table[x={N}, y={err_kryl}] {Experiments/experiment4/errors_final.dat};
\addlegendentry{\small \texttt{IEP\_KRYL}}
\addplot[green!50!black, mark=triangle*, solid] table[x={N}, y={err_kryl_reorth}] {Experiments/experiment4/errors_final.dat};
\addlegendentry{\small  \texttt{IEP\_KRYLREORTH(full)}}

\nextgroupplot[label style={font=\tiny}, tick label style={font=\tiny}, mark options={scale=0.6},legend columns = 1,legend style={at = {(0.01,0.99)}, anchor = north west,row sep=-0.2pt,nodes={scale=0.55, transform shape}},legend cell align={left},grid=major, xlabel={$N$},ylabel={$\|W^TV-I_N\|_2$}, xmin = 5, xmax = 50, ymin=1e-17,ymax=1e11, ymode=log, xtick distance=10, ytick distance = 10000,minor ytick={1e-16, 1e-15, 1e-14, 1e-12, 1e-11, 1e-10, 1e-8, 1e-7, 1e-6, 1e-4, 1e-3, 1e-2, 1e0, 1e1, 1e2, 1e4, 1e5, 1e6, 1e8, 1e9, 1e10, 1e12, 1e13, 1e14, 1e16, 1e17, 1e18, 1e20, 1e21, 1e22, 1e24, 1e25, 1e26,1e28, 1e29, 1e30,1e32, 1e33, 1e34}]

\addplot[red, mark=*, solid] table[x={N}, y={err_bior_core}] {Experiments/experiment4/errors_final.dat};
\addlegendentry{\small  \texttt{IEP\_CORE}}
\addplot[blue, mark=diamond*, solid] table[x={N}, y={err_bior_kryl}] {Experiments/experiment4/errors_final.dat};
\addlegendentry{\small \texttt{IEP\_KRYL}}
\addplot[green!50!black, mark=triangle*, solid] table[x={N}, y={err_bior_kryl_reorth}] {Experiments/experiment4/errors_final.dat};
\addlegendentry{\small  \texttt{IEP\_KRYLREORTH(full)}}
\addplot[red, mark=diamond*, dashed] table[x={N}, y={err_bior_core_scaled}] {Experiments/experiment4/errors_final.dat};
\addlegendentry{\small \texttt{IEP\_CORE} \textbf{scaled}}
\addplot[green!50!black, mark=triangle*, dashed] table[x={N}, y={err_bior_kryl_reorth_scaled}] {Experiments/experiment4/errors_final.dat};
\addlegendentry{\small  \texttt{IEP\_KRYLREORTH(full) \textbf{scaled}}}

\end{groupplot}
\end{tikzpicture}
    \caption{Error comparison (left) and loss of biorthogonality (right) for Chebyshev nodes in $[-1,1]$ and random weights, averaged over 20 runs.}
    \label{fig:experiment4}
\end{figure}
\end{itemize}

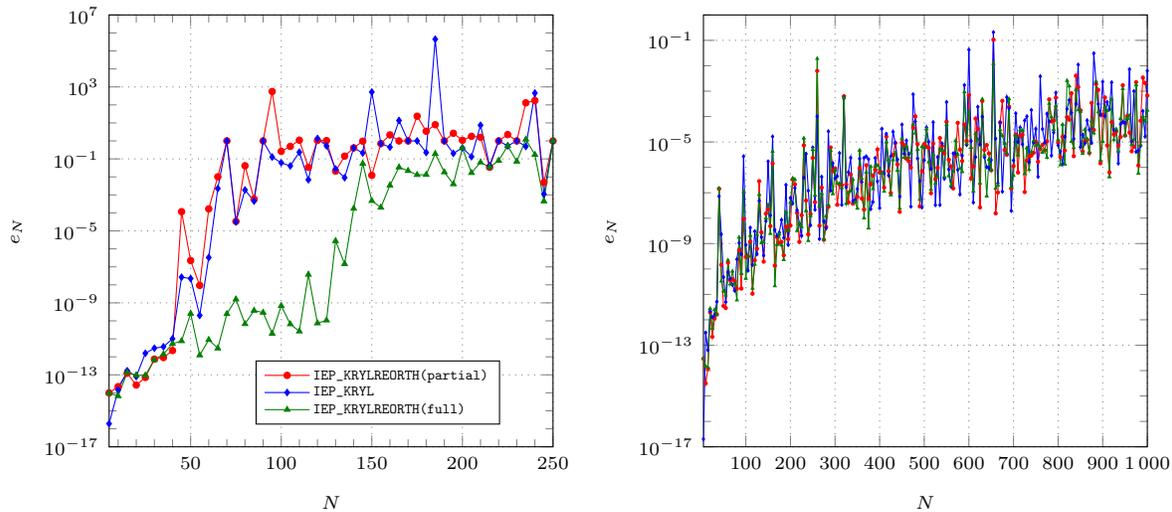
\begin{figure}
    \centering
    \pgfkeys{/pgf/number format/1000 sep={\,}}
    \tikzsetnextfilename{experiment5}
\pgfplotsset{height=0.47\linewidth,width=0.48\linewidth}
\pgfplotsset{major grid style={dotted,gray}}

\noindent%
\begin{tikzpicture}[scale=1,baseline]%
\begin{groupplot}[
    group style={
        group size=2 by 1,
        horizontal sep=2cm,
    }
]

\nextgroupplot[label style={font=\tiny}, tick label style={font=\tiny},mark options={scale=0.6}, legend columns = 1, legend style={at = {(0.33,0.2)}, anchor = north west,row sep=-0.2pt,nodes={scale=0.55, transform shape}},legend cell align={left}, grid=major, xlabel={$N$},ylabel={$e_N$}, xmin = 5, xmax = 250, ymin=1e-17,ymax=1e7, ymode=log, xtick distance=50, minor x tick num = 4, ytick distance = 10000, minor ytick={1e-16, 1e-15, 1e-14, 1e-12, 1e-11, 1e-10, 1e-8, 1e-7, 1e-6, 1e-4, 1e-3, 1e-2, 1e0, 1e1, 1e2, 1e4, 1e5, 1e6, 1e8, 1e9, 1e10, 1e12, 1e13, 1e14, 1e16, 1e17, 1e18, 1e20, 1e21, 1e22, 1e24, 1e25, 1e26,1e28, 1e29, 1e30,1e32, 1e33, 1e34}]

\addplot[red, mark=*, solid] table[x={N}, y={err_krylpartial}] {Experiments/experiment5/data_experiment5_equid2.dat};
\addlegendentry{\small \texttt{IEP\_KRYLREORTH(partial)}}
\addplot[blue, mark=diamond*, solid] table[x={N}, y={err_krylnoreorth}] {Experiments/experiment5/data_experiment5_equid2.dat};
\addlegendentry{\small \texttt{IEP\_KRYL}}
\addplot[green!50!black, mark=triangle*, solid] table[x={N}, y={err_krylfull}] {Experiments/experiment5/data_experiment5_equid2.dat};
\addlegendentry{\small  \texttt{IEP\_KRYLREORTH(full)}}

\nextgroupplot[label style={font=\tiny}, tick label style={font=\tiny}, mark options={scale=0.45},legend columns = 1,legend style={at = {(0.4,0.2)}, anchor = north west,row sep=-0.2pt,nodes={scale=0.55, transform shape}},legend cell align={left},grid=major, xlabel={$N$},ylabel={$e_N$}, xmin = 5, xmax = 1000, ymin=1e-17,ymax=1e0, ymode=log, xtick distance=100, ytick distance = 10000,minor ytick={1e-16, 1e-15, 1e-14, 1e-12, 1e-11, 1e-10, 1e-8, 1e-7, 1e-6, 1e-4, 1e-3, 1e-2, 1e0, 1e1, 1e2, 1e4, 1e5, 1e6, 1e8, 1e9, 1e10, 1e12, 1e13, 1e14, 1e16, 1e17, 1e18, 1e20, 1e21, 1e22, 1e24, 1e25, 1e26,1e28, 1e29, 1e30,1e32, 1e33, 1e34}]

\addplot[red, mark=*, solid] table[x={N}, y={err_krylpartial}] {Experiments/experiment5/data_experiment5_cheb2.dat};
\addplot[blue, mark=diamond*, solid] table[x={N}, y={err_krylnoreorth}] {Experiments/experiment5/data_experiment5_cheb2.dat};
\addplot[green!50!black, mark=triangle*, solid] table[x={N}, y={err_krylfull}] {Experiments/experiment5/data_experiment5_cheb2.dat};

\end{groupplot}
\end{tikzpicture}
    \caption{Error comparison for equidistant nodes in [-1,1] and uniform random weights (left) and Chebyshev nodes in [-1,1] and uniform random weights (right) for larger values of $N$}
    \label{fig:experiment5}
\end{figure}

In a final experiment, illustrated in \cref{fig:experiment5}, we compare the precision of the Krylov-based algorithms for larger values of $N$, \rev{even though $N$ is typically only moderate in practical applications~\cite{LauMastVanAVanDooren24}}. More precisely, \texttt{IEP\_KRYL},\texttt{IEP\_KRYLREORTH(full)} and \texttt{IEP\_KRYLREORTH(partial)} are compared using the nodes and weights described above. For the well-conditioned Chebyshev example, all three methods have comparable performance, with a relative forward error $e_N$ of approximately $10^{-4}$ at $N = 1000$. In contrast, for the case with equidistant nodes in the interval $[-1,1]$, the benefits of full reorthogonalization become apparent. Note that, up to $N=50$, \texttt{IEP\_KRYL} and \texttt{IEP\_KRYLREORTH(partial)} have acceptable accuracy and that the fully reorthogonalized variant \texttt{IEP\_KRYLREORTH(full)} extends this limit to $N = 150$ for a similar level of relative error.

\section{Conclusion}\label{sec:conclusion}
By connecting MOPs, Krylov subspaces and inverse eigenvalue problems, we derive a Krylov-based method and a core transformation algorithm to compute the unique recurrence coefficients associated with the MOPs on the step-line. The methods \texttt{IEP\_CORE} and \texttt{IEP\_KRYLREORTH(full)} have similar performance for the ill-conditioned examples associated with Hahn and Kravchuk. However, the experiments on better-conditioned problems show that \texttt{IEP\_KRYLREORTH(full)} achieves the highest accuracy and overall best performance among the methods considered, at the price of a higher computational cost.  A theoretical analysis of the numerical properties of the IEP could provide more insight, for example,
on the conditioning of the problem.
Furthermore, an
updating/downdating approach can be applied as well,
which has proven to be 
more flexible as it allows reusing the existing recurrence relations to build new recurrences under modest modifications, such as adding or deleting nodes of
the inner product \cite{GraHar84,Rut63,Vanb23,VanBVanBVand24,MR3867618,VanbVanbVand22}.
\rev{Another possibility, that we did not discuss, is to extend Laurie's approach~\cite{Laurie99} for classical orthogonal polynomials based on continued fractions. For MOPs, one can look into the connections with vector continued fractions~\cite{Kaliaguine95}, or the recent developments on branched continued fractions~\cite{Sokal24}.}
A natural direction for future research is the generalization towards other recurrence relations for MOPs such as the nearest-neighbor recurrence relations \cite{Van-nearestneighbor}. Finally, it would be interesting to compare 
our methods with the techniques in~\cite{MilStan03,FilHanVanass15}. 
\section*{Funding}
The research was partially supported by the Research Council KU Leuven (Belgium), project C16/21/002 (Manifactor: Factor Analysis for Maps into Manifolds) and by the Fund for Scientific Research -- Flanders (Belgium), projects G0A9923N (Low rank tensor approximation techniques for up- and downdating of massive online time series clustering) and G0B0123N (Short recurrence relations for rational Krylov and orthogonal rational functions inspired by modified moments), and junior postdoctoral fellowship 12A1325N (Short Recurrences for Block Krylov Methods with Applications to Matrix Functions and Model Order Reduction) for the second author.

\rev{\section*{Acknowledgements} We would like to thank Walter Van Assche for his insightful suggestions on this topic. We also thank the referees for their valuable comments, which helped improve the quality of this manuscript.}

\bibliographystyle{siamplain}
\bibliography{references}
\end{document}